\documentclass[11pt]{article}
\usepackage[T1]{fontenc}
\usepackage[utf8]{inputenc}
\usepackage{lmodern}
\usepackage[a4paper,margin=31mm]{geometry}
\usepackage{amsmath,amssymb,amsthm,mathtools,mathrsfs}
\usepackage{booktabs,array,graphicx}
\usepackage{enumitem}
\usepackage{flafter}
\usepackage[section]{placeins}
\usepackage{microtype}
\usepackage{xcolor}
\usepackage{needspace}
\usepackage{hyperref}
\hypersetup{
  colorlinks=true,
  linkcolor=blue!55!black,
  citecolor=blue!55!black,
  urlcolor=blue!55!black,
  pdftitle={An Odd-Arch Family of Perturbed Hofstadter Q-Recursions: Gaussian Profiles and Exact Frequency Laws},
  pdfauthor={Marco Mantovanelli},
  pdfsubject={Odd-arch family, Gaussian profiles, finite transitions, parameter-dependent frequencies and exact sums},
  pdfkeywords={Hofstadter recursion, meta-Fibonacci sequence, binary word,
    Gaussian isoperimetric profile, dyadic frequency law}
}
\setlist{itemsep=2pt,topsep=5pt}
\numberwithin{equation}{section}
\newtheorem{theorem}{Theorem}[section]
\newtheorem{lemma}[theorem]{Lemma}
\newtheorem{proposition}[theorem]{Proposition}
\newtheorem{corollary}[theorem]{Corollary}
\theoremstyle{definition}
\newtheorem{definition}[theorem]{Definition}
\theoremstyle{remark}
\newtheorem{remark}[theorem]{Remark}

\newcommand{\I}{\mathcal I}
\newcommand{\Sp}{\mathcal S_+}
\newcommand{\Sm}{\mathcal S_-}
\newcommand{\eps}{\varepsilon}

\title{An Odd-Arch Family of Perturbed Hofstadter $Q$-Recursions:\\
Gaussian Profiles and Exact Frequency Laws}
\author{Marco Mantovanelli\\[.35em]
\small\href{mailto:Marco@Mantovanelli.de}{\texttt{Marco@Mantovanelli.de}}\\[-.1em]
\small\href{https://orcid.org/0009-0002-0631-293X}
{ORCID: 0009-0002-0631-293X}}
\date{}
\begin{document}
\maketitle
\begin{abstract}
For every integer $p\ge 2$, we construct a globally defined
positive odd sequence $Q_p$ from an explicit finite initial
block, satisfying the parity-perturbed Hofstadter recursion
\[
Q_p(n)
=Q_p\bigl(n-Q_p(n-1)\bigr)
+Q_p\bigl(n-Q_p(n-2)\bigr)-(-1)^n.
\]

A causal binary-word construction yields the exact factorization
\[
U_{p,k}
=\bigl(A_k\overline{A}_{k+1}\bigr)^{p-1}A_{k+2},
\]
producing exactly $2p-1$ alternating signed excursions in
each natural block.

The normalized elementary heights converge uniformly to the
Gaussian isoperimetric profile
$J=\varphi\circ\Phi^{-1}$.
Consequently, every fixed $p$ admits an explicit dyadic
limiting profile with $2p-1$ sign transitions and a sharp
relative-error envelope:
\[
\begin{aligned}
\limsup_{n\to\infty}
\sqrt{\log_2 n}
\left|\frac{Q_p(n)}{n}-\frac12\right|
&=\Lambda_p=\frac{b_p}{2},\\
b_p\left(\frac{3p-1}{2}+\Phi(-b_p)\right)
&=\varphi(b_p),
\end{aligned}
\]
where $\varphi$ and $\Phi$ denote the standard normal density
and distribution function.

A refined low-height analysis determines all zeros and
parity crossings in every sufficiently large dyadic block.
For $p\ge3$, such blocks contain exactly $2p-1$ transverse
even-index zeros, whereas the three-arch case $p=2$
exhibits an alternating four/six-zero pattern.

An ordered factorization of value frequencies yields
explicit dyadic histograms, exact moments, total-variation
distances, and the positions of extremal multiplicities.
The limiting frequency generating function
$G(z)=z^3/(2-z)$ is universal, although the highest
finite multiplicities depend on $p$.

Finally, writing $C_p=2(3p+1)/3$, we establish the
summatory asymptotic
\[
\begin{aligned}
\sum_{n\le N}Q_p(n)
={}&\frac{N^2}{4}
+\frac{N}{2}\log_2 N\\
&+N\Psi_p\left(
\left\{\log_2\frac{N}{C_p}\right\}
\right)
+O_p\left(\frac{N}{\sqrt{\log N}}\right),
\end{aligned}
\]
where $\Psi_p$ is an explicitly determined continuous
one-periodic function.
\end{abstract}



\smallskip
\noindent\textbf{Supplementary material.}
Reproducibility archive:
\href{https://doi.org/10.5281/zenodo.23593788}{10.5281/zenodo.23593788}.

\section{Introduction and main results}\label{sec:intro}
Changing a finite initial block in a nested recurrence changes which
parts of the past are read, not merely the first few values. Here this
mechanism produces an explicit family $Q_p$, indexed by integers $p\ge2$,
with any prescribed odd number $2p-1$ of signed excursions per natural
scale. All members use the same recurrence; the parameter enters only
through the initial block. Their elementary excursions have a common
Gaussian limit, while their arrangement, finite transitions, and upper
value multiplicities depend on the parameter.

\subsection{The family and its global structure}
Write
\[
 n_p=\begin{cases}8,&p=2,\\6p-6,&p\ge3.\end{cases}
\]
We prescribe every value $Q_p(n)$ for $1\le n<n_p$ and then apply
\begin{equation}\label{eq:family-rec}
 Q_p(n)=Q_p\bigl(n-Q_p(n-1)\bigr)
       +Q_p\bigl(n-Q_p(n-2)\bigr)-(-1)^n,
 \qquad n\ge n_p.
\end{equation}
For $p=2$ the prescribed block is $(1,1,3,3,3,3,5)$.
The three-arch member $Q_2$ is recorded in the On-Line
Encyclopedia of Integer Sequences as
\href{https://oeis.org/A401334}{A401334} \cite{oeisA401334}.
For $p\ge3$, put $m=\lceil n/2\rceil$ and prescribe, for
$1\le n\le6p-7$,
\begin{equation}\label{eq:family-seed}
 Q_p(n)=
 \begin{cases}
 2\min(m-1,p-1)+1,&m<p+2,\\
 2p+2\lceil(m-p-2)/2\rceil+1,&m\ge p+2,\ n\text{ odd},\\
 2p+2\lfloor(m-p-2)/2\rfloor+1,&m\ge p+2,\ n\text{ even}.
 \end{cases}
\end{equation}
The entire block is part of the definition: the recurrence must not be
applied before $n_p$. Each recursive argument must belong to
$\{1,\ldots,n-1\}$; no values at zero or at negative indices are supplied.

\begin{theorem}[Global odd-arch family]\label{thm:family-global}\label{thm:global}
For every integer $p\ge2$, the initial block above and
\eqref{eq:family-rec} determine a unique sequence for all positive indices.
All values are odd, and
\begin{equation}\label{eq:family-basic}
 1\le Q_p(n)\le n,\qquad Q_p(n+2)-Q_p(n)\in\{0,2\}
 \qquad(n\ge1).
\end{equation}
Each parity subsequence assumes every positive odd value, each only
finitely often. The independently constructed block words satisfy
\begin{equation}\label{eq:family-factor}
 U_{p,k}=(A_k\overline{A}_{k+1})^{p-1}A_{k+2}\qquad(k\ge2),
\end{equation}
where the primitive Dyck words $A_k$ are defined in
\eqref{eq:Awords} and have length
$\ell_k=(2^k+2(-1)^k)/3$.
A bar is bitwise complementation, not reversal. Consequently each
natural block has exactly $2p-1$ primitive signed excursions, with
asymptotic relative lengths $1,2,\ldots,1,2,4$.
\end{theorem}

The proof in Section~\ref{sec:family} constructs the words and two slow
cores before checking any nested identity. In particular, existence is
not assumed in order to establish the block structure. The first three
members are the three-, five-, and seven-arch sequences $Q_2,Q_3,Q_4$:
\begin{center}
\small
\begin{tabular}{@{}clc@{}}\toprule
Member & Prescribed initial block & Recursion starts at\\\midrule
$Q_2$ & $(1,1,3,3,3,3,5)$ & $8$\\
$Q_3$ & $(1,1,3,3,5,5,5,5,7,7,9)$ & $12$\\
$Q_4$ & $(1,1,3,3,5,5,7,7,7,7,9,9,11,9,11,11,13)$ & $18$\\\bottomrule
\end{tabular}
\end{center}

\subsection{Shape, finite transitions, frequencies, and summation}
The common elementary profile is
\[
 J(x)=\varphi(\Phi^{-1}(x)),\qquad 0<x<1,\qquad J(0)=J(1)=0,
\]
where $\varphi$ and $\Phi$ are the standard normal density and distribution
function. Section~\ref{sec:gaussian} proves uniform convergence of the
normalized elementary heights. Section~\ref{sec:family-profiles} combines
these limits at the exact block positions and obtains an explicit
standard dyadic profile $\mathcal P_p$. For each fixed $p$ and index
parity $r$,
\[
 (-1)^{j+r}\sqrt j\left(\frac{Q_p(n)}n-\frac12\right)
 \longrightarrow \mathcal P_p(n/2^j)
\]
uniformly on $2^j\le n\le2^{j+1}$ along that parity, in the precise
sampling convention of Theorem~\ref{thm:family-profile}.
There are exactly $2p-1$ limiting sign transitions. The sharp envelope is
\[
 \limsup_{n\to\infty}\sqrt{\log_2n}
 \left|\frac{Q_p(n)}n-\frac12\right|=\Lambda_p=\frac{b_p}{2},\qquad
 b_p\left(\frac{3p-1}{2}+\Phi(-b_p)\right)=\varphi(b_p),
\]
with a unique positive root $b_p$. In particular,
$\Lambda_p\sim1/(3p\sqrt{2\pi})$ as $p\to\infty$.
This last assertion concerns the constants after the fixed-$p$ limit;
no simultaneous limit in $p$ and $n$ is asserted.

The finite structure contains more information than the profile alone.
Section~\ref{sec:family-boundaries} gives exact finite zero and crossing
tables. For $p\ge3$ and $j\ge p+\lfloor\log_2p\rfloor+10$, the $j$-th
standard block has exactly $2p-1$ even-index zeros, all transverse.
The four/six-zero alternation at $p=2$ is a genuine exceptional feature,
not a universal rule.

For $\nu_p(v)=\#\{n\ge1:Q_p(n)=v\}$, define
\[
 D_{p,j}(z)=\sum_{u=2^j}^{2^{j+1}-1}z^{\nu_p(2u+1)}.
\]
Section~\ref{sec:family-frequencies} derives both the ordered frequency
words and their histograms. The normalized limit is always $z^3/(2-z)$,
but the finite upper tail changes with $p$. With
$d=\lfloor\log_2p\rfloor$, the exact moments on $j\ge d+2$ are
\[
 \mathbb E X_{p,j}=4-2^{-j},\qquad
 \operatorname{Var}(X_{p,j})
 =2-\bigl(2(j-d)+4p-5\bigr)2^{-j}-4^{-j}.
\]
Theorem~\ref{thm:family-dyadic} gives all coefficients, split according
to whether $p$ is a power of two. The following results determine exact
total-variation distances and every location of the two highest
multiplicities. A worked $p=3$ example makes the join correction and the dyadic
cuts explicit before the general histogram formulas.

Finally, Section~\ref{sec:family-summation} proves an exact telescoping
identity and, with $C_p=2(3p+1)/3$,
\[
 \sum_{n\le N}Q_p(n)=\frac{N^2}{4}+\frac N2\log_2N
 +N\Psi_p(\{\log_2(N/C_p)\})+O_p(N/\sqrt{\log N}),
\]
where
$\Psi_p(\theta)=p-1-\tfrac12\log_2C_p-\tfrac\theta2-2^{-\theta}$
on $[0,1]$, extended continuously with period one.
Thus the logarithmic bias and the periodic shape are shared by the
family, while the phase scale and vertical constant are explicit
functions of $p$.

\subsection{Context and organization}
Hofstadter introduced the unperturbed $Q$-recursion in
\emph{G\"odel, Escher, Bach} \cite{hofstadter1979}.
For related nested recurrences, explicit combinatorial descriptions
have often supplied the missing control of the self-referential
arguments. Tanny \cite{tanny1992} analyzed a monotone relative;
Balamohan, Kuznetsov, and Tanny \cite{balamohan2007} obtained a complete
description of another variant through its generations and value
frequencies. Tree representations provide a further source of exact
solutions \cite{jacksonruskey2006}.
These examples suggest looking for a finite combinatorial evolution
behind the numerical sequence.

The original parity-perturbed recursion was introduced by the present
author \cite{mantovanelli} and is recorded as OEIS A394051 \cite{oeis}.
It starts with $(1,1)$ and adds $+(-1)^n$ to the two nested terms.
Cloitre \cite{cloitre} proved its global definedness and the optimal
error order $n/\sqrt{\log n}$ by an exact-fit interleaving construction,
binary heights, and canonical plane forests. The same framework led
to a dyadic frequency law \cite{mantovanelli}. A subsequent analysis
of the associated two-variable zeta function identified a periodic
summatory term and a weak Gaussian limit for a distinguished family
of layers \cite{mantovanellizeta}.

A related delayed perturbation admits uniform Gaussian profiles through
nonlinear averaging and smooth comparison functions
\cite{mantovanellidelayed}. Here the elementary lengths alternate between
$2L+2$ and $2L-2$. An exact change of coordinates removes both the length
correction and the additive height term, allowing that comparison
method to treat the two maps within one identity.

Our construction uses Cloitre's exact-fit framework. In particular, the
tape counters and height additivity come from
\cite[Definition 2.3 and Lemmas 2.4--2.5]{cloitre}, and the parked-head
verification follows \cite[Lemmas 3.3--3.4]{cloitre}. On Dyck inputs,
$\Sp$ agrees with Cloitre's transform $T$; if $P=0R1$ is primitive Dyck,
then $\Sm(P)=T(R)$ \cite[Section 4.1]{cloitre}. The elementary excursions
are inherited. The contribution here is their repeated signed-pair
organization, its explicit finite seeds and block positions, and the
resulting parameter-dependent shape, crossing, frequency, and summation
laws. Section~\ref{sec:comparison} gives the exact identification with
the original word family and separates the inherited tools from these
results.

The main text develops the general parameter $p$ throughout.
Appendix~\ref{app:three-arch} records the detailed $p=2$ formulas,
including its exceptional tangencies, pure valuation histogram, and
worked finite crossing example. This retains the three-arch model as
an explicit specialization rather than as a prerequisite for the
familial construction.
\begin{figure}[htbp]
\centering
\includegraphics[width=.99\textwidth]{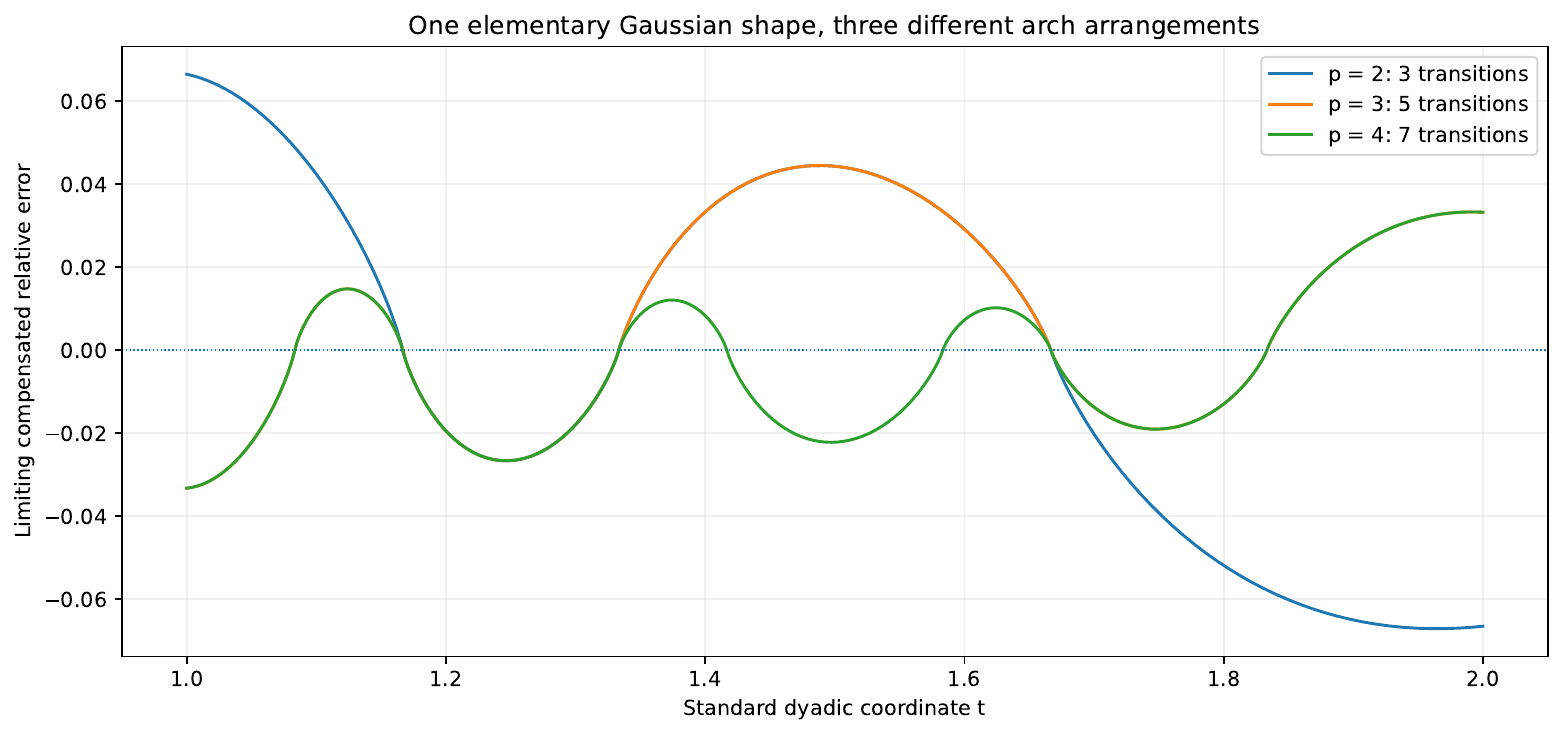}
\caption{The three explicit limiting standard-block profiles for
$p=2,3,4$. The same elementary Gaussian function produces distinct
arrangements and different sharp amplitudes.}
\label{fig:family-all-profiles}
\end{figure}

\section{Causal word maps and elementary excursions}\label{sec:words}

For a binary word $P$ of length $L$, let $o_P(t)$ and $z_P(t)$ count
ones and zeros in its prefix of length $t$, and let
\[
h_P(t)=z_P(t)-o_P(t),\qquad 0\le t\le L.
\]
Positions in words are indexed from zero. A word is balanced when it has
equally many zeros and ones. It is a \emph{primitive Dyck word} when it
is balanced and $h_P(t)>0$ at every integer $0<t<L$.
A bar denotes bitwise complementation, not reversal.

The following strict machine is the exact-fit case of
\cite[Definition 2.3]{cloitre}. No exhausted-tape fallback is needed
for the words used here.

\begin{definition}[Strict interleaving and the two maps]
The machine $\I_0(X,Y)$ starts in state zero. In state zero it reads
and emits the next symbol of $X$; in state one it does the same with
$Y$. The emitted bit is the next state. The machine stops when both
tapes have been consumed. Reading a designated tape which is already
empty before that time is a failure: there is no fallback rule.
For a balanced word beginning in zero and ending in one, define
\begin{equation}\label{eq:maps}
\Sp(P)=\I_0(0P1,P),\qquad
\Sm(P)=\I_0(P,P[1:-1]).
\end{equation}
Here $P[1:-1]$ deletes its initial zero and final one.
The first map also applies to the balanced seeds
$(10)^{p-1}01$ for $p\ge2$; the exceptional seed $1001$ is the case $p=2$.
\end{definition}

\begin{lemma}[Causality and primitive excursions]\label{lem:maps}
Both maps in \eqref{eq:maps} are defined, balanced, begin in zero and
end in one. Their lengths are $2L+2$ and $2L-2$, respectively.
Both preserve primitive Dyck words. The map $\Sp$ is also defined on
any balanced word ending in one, regardless of its first symbol.
\end{lemma}
\begin{proof}
After a nonempty output prefix let $i,j$ be the tape counters, $z,o$
the output counts, and $c$ the last bit. Counting the states before
successive reads, as in \cite[Lemma 2.5]{cloitre}, gives
\begin{equation}\label{eq:counters}
i=z+c,\qquad j=o-c.
\end{equation}
For $\Sp$, the total number of zeros and of ones is $L+1$.
A request to read an exhausted zero tape would have $i=L+2,c=0$,
forcing too many zeros. An exhausted one tape with $j=L,c=1$ forces
$o=L+1$: all ones have been consumed. Since the zero tape ends in one,
that tape must then also be exhausted. This is the legitimate stopping
state. For $\Sm$, replace these tape lengths by $L,L-2$ and the total
counts by $L-1$; the same argument applies. Equation~\eqref{eq:counters}
at completion forces the final bit to be one. The first bit is zero.

For $\Sp$ the prefix height, while $1\le i\le L+1$, is
$1+h_P(i-1)+h_P(j)>0$ if $P$ is Dyck. If the zero tape has already
finished, all zeros have been emitted, and the height equals the
number of remaining output symbols. It is positive until completion.
For $\Sm$ the height is
\[
h_P(i)+h_P(j+1)-1.
\]
Before the zero tape finishes, $1\le i<L$ and $1\le j+1\le L-1$;
primitivity makes both parent heights at least one. After it finishes,
the same remaining-length argument applies. This proves primitivity.
\end{proof}

Let $\alpha_i$ be the number of zeros immediately preceding the $i$-th
one, back to the previous one or the word's beginning. Let $\beta_i$
be the number of ones after the $i$-th zero, up to the next zero or the
end. The gaps may be zero. A balanced word of length $2p$ has $p$ zeros
and therefore $p$ gaps of each kind.

\Needspace{14\baselineskip}
\begin{lemma}[Run expansions and symmetry]\label{lem:runs}
Let $P$ be a balanced binary word of length $2p$, beginning in zero
and ending in one. For $p\ge1$,
\begin{equation}\label{eq:plus-runs}
\Sp(P)=0^{\alpha_1+1}1
\left(\prod_{i=1}^{p-1}0^{\alpha_{i+1}+1}1^{\beta_i+1}\right)
01^{\beta_p+1}.
\end{equation}
For $p\ge2$,
\begin{equation}\label{eq:minus-runs}
\Sm(P)=0^{\alpha_1}1^{\beta_1+1}
\left(\prod_{i=2}^{p-1}0^{\alpha_i+1}1^{\beta_i+1}\right)
0^{\alpha_p+1}1^{\beta_p}.
\end{equation}
For $P=01$, one has $\Sm(P)=01$.
If $P$ equals its reverse-complement, both outputs do too.
The initial zero-run and terminal one-run both gain one symbol under
$\Sp$ and retain their lengths under $\Sm$.
\end{lemma}
\begin{proof}
This is the two-map form of the gap transport in \cite[Lemma 4.1]{cloitre}.
A zero emitted from the one tape starts a new zero-run, and the zero
tape then emits its next zero gap followed by the one ending that gap.
The one tape continues the corresponding one-run. Tracking the two
initializations and the tape ends gives the displayed expansions.
Reverse-complement symmetry of $P$ is equivalent to
$\alpha_i=\beta_{p+1-i}$. The displayed run lengths then match in reverse
order with zeros and ones exchanged. The outer-run statements follow
from the first and last factors.
\end{proof}

Define the elementary words and their lengths by
\begin{equation}\label{eq:Awords}
A_2=01,\qquad
A_{k+1}=\begin{cases}\Sp(A_k),&k\text{ odd},\\
\Sm(A_k),&k\text{ even},\end{cases}\qquad
\ell_k=|A_k|=\frac{2^k+2(-1)^k}{3}\quad(k\ge2).
\end{equation}
Thus $A_3=01$, $A_4=001011$, and $A_5=0010011011$.
The length formula follows from $\ell_{k+1}=2\ell_k-2(-1)^k$.
All $A_k$ are primitive and reverse-complement symmetric. Their initial
zero-run and terminal one-run have length $\lfloor k/2\rfloor$.

\begin{lemma}[Concatenation transport]\label{lem:factor-transport}
For balanced words $P,R,S$ beginning in zero and ending in one,
\begin{align}
\Sp(P\overline R S)&=\Sp(P)\,\overline{\Sm(R)}\,\Sp(S),\label{eq:factor-plus}\\
\Sm(P\overline R S)&=\Sm(P)\,\overline{\Sp(R)}\,\Sm(S).\label{eq:factor-minus}
\end{align}
\end{lemma}
\begin{proof}
In \eqref{eq:factor-plus}, the first output segment uses $0P1$ from the
zero tape and $P$ from the one tape. The terminal one on the first tape
is borrowed from the first bit of $\overline R$.
The state is now one. Complementing the output and exchanging the
roles of the tapes makes the next segment exactly $\Sm(R)$: the two
consumed pieces are $R$ and its interior. This leaves the last zero of
$\overline R$ on the original zero tape. The state returns to zero and
the remaining tapes are $0S1,S$, giving $\Sp(S)$.

For \eqref{eq:factor-minus}, the first segment consumes $P$ and the
interior of $P$, leaving its terminal one on the one tape. In state
one, complementation and tape exchange give $\Sp(R)$; its extra terminal
one is borrowed from the complement of the initial zero of $S$.
The remaining original tapes are $S,S[1:-1]$ in state zero.
All segment computations are causal by Lemma~\ref{lem:maps}, so these
concatenations also verify the strict reading convention at each join.
\end{proof}

\section{Construction and global existence of the family}\label{sec:family}\label{sec:cores}
\subsection{Block positions and independent cores}
Fix an integer $p\ge2$ and set
\begin{equation}\label{eq:family-skeleton}
\begin{split}
 c_p&=3p+1,\qquad C_p=\frac{2c_p}{3},\qquad
 U_{p,1}=(10)^{p-1}01,\\
 L_{p,k}&=\frac{c_p2^k+2(-1)^k}{3},\qquad
 a_{p,k}=\frac{c_p2^k-3-(-1)^k}{6},\\
 d_{p,k}&=p-2+\left\lfloor\frac{k-1}{2}\right\rfloor,
 \qquad s_{p,k}=2a_{p,k}-d_{p,k}.
\end{split}
\end{equation}
Advance $U_{p,k}$ by $\Sp$ for odd $k$ and by $\Sm$ for even $k$.
The two cores are first defined independently of any nested recursion:
\[
 F_p(m)=G_p(m)=\min(m-1,p-1)\qquad(1\le m<p+2).
\]
On the $k$-th block, with $o=o_{U_{p,k}}$ and $z=z_{U_{p,k}}$, put
\begin{align}
 F_p(s_{p,k}+t)&=a_{p,k}+o(t),&
 G_p(s_{p,k}+t)&=a_{p,k}+z(t)
 &&\left\{\begin{array}{l}k\text{ odd},\\0\le t\le L_{p,k},\end{array}\right.
 \label{eq:family-core-odd}\\
 F_p(s_{p,k}+t)&=a_{p,k}+z(t),&
 G_p(s_{p,k}+t)&=a_{p,k}+o(t+1)
 &&\left\{\begin{array}{l}k\text{ even},\\0\le t<L_{p,k}.\end{array}\right.
 \label{eq:family-core-even}
\end{align}
Finally define
\begin{equation}\label{eq:family-lift}
 Q_p(2m-1)=2F_p(m)+1,\qquad Q_p(2m)=2G_p(m)+1.
\end{equation}
The case $p=2$ is the three-arch sequence detailed in
Appendix~\ref{app:three-arch}.

The initial word is balanced and ends in one. Its first plus transform
is causal by Lemma~\ref{lem:maps}; all later words begin in zero and
end in one, so both maps remain causal. The definitions satisfy
\begin{equation}\label{eq:family-clock-identities}
\begin{split}
 a_{p,k+1}&=a_{p,k}+L_{p,k}/2,\\
 s_{p,k+1}-s_{p,k}&=L_{p,k}-\mathbf1_{\{k\text{ even}\}},\\
 s_{p,k}-a_{p,k}&=s_{p,k-1}\qquad(k\ge2).
\end{split}
\end{equation}
Consequently the two core definitions agree at shared endpoints.
At $s_{p,1}=p+2$ both cores rise from $p-1$ to $p$. All subsequent
increments are zero or one. The cores are unbounded and satisfy
\begin{equation}\label{eq:core-bound}
\begin{gathered}
0\le F_p(m),G_p(m)\le m-1\quad(m\ge1),\\
1\le F_p(m),G_p(m)\le m-2\quad(m\ge p+1).
\end{gathered}
\end{equation}
The second upper bound follows from the repeated value $p-1$ at
$m=p,p+1$ and the unit-step property. No nested recursion has been used.

\subsection{Exact coupled recursions}
\begin{proposition}[Coupled identities on every subsequent block]\label{thm:cores}
For every $p\ge2$ and $m>s_{p,2}=3p+2$, the word-defined cores satisfy
\begin{align}
 F_p(m)&=G_p(m-G_p(m-1)-1)+G_p(m-F_p(m-1)-1)+1,
 \label{eq:Frec}\\
 G_p(m)&=F_p(m-F_p(m))+F_p(m-G_p(m-1)).\label{eq:Grec}
\end{align}
All arguments lie in $\{1,\ldots,m-1\}$.
\end{proposition}
\begin{proof}
Fix $p$ and suppress it temporarily in $F,G,U_k,a_k,s_k,L_k$.
The bounds in \eqref{eq:core-bound} put the two arguments in
\eqref{eq:Frec} in $[2,m-2]$ and those in \eqref{eq:Grec} in
$[2,m-1]$ and $[3,m-1]$, respectively. We verify the identities on
$s_k<m\le s_{k+1}$ for every $k\ge2$, assigning a shared endpoint
to its left block. Both core formulas equal $a_{k+1}$ at that endpoint,
so this convention agrees with their independent definition.

Here is the precise parked-head observation used below, in the form
of \cite[Lemmas 3.3--3.4]{cloitre}. Suppose that $o_-,z_-$ are the
counts before the $t$-th output, with $t\ge1$. The counters immediately
after that output are
\begin{equation}\label{eq:read-counters}
i=z_-+1,\qquad j=o_-.
\end{equation}
Indeed, \eqref{eq:counters} gives the counters before the read; the
read increments $i$ if the old state is zero and $j$ if it is one.
This also holds for $t=1$, when both counters and the old state are
zero. If the new output bit is one, the last symbol consumed from the
zero tape is one: it either produced this output or initiated the
current succession of reads from the one tape. During those reads
the zero tape remains parked at that symbol.

In the following table, $s'=s_{k-1}$, and the entries are the inner
arguments minus $s'$. We abbreviate the child prefix counts by
$o=o(t)$, $z=z(t)$, $o_-=o(t-1)$, and $z_-=z(t-1)$.
\begin{center}
\begin{tabular}{@{}cccc@{}}\toprule
Child level & Map & $F$-identity arguments & $G$-identity arguments\\\midrule
odd & $\Sm$ & $(o_-,z_-)$ & $(z,o_-+1)$\\
even & $\Sp$ & $(z-1,o_-)$ & $(o,z)$\\\bottomrule
\end{tabular}
\end{center}
We verify both the values and the domains in each row.

\smallskip\noindent\emph{Odd $k\ge3$.}
Let $P=U_{k-1}$, $\ell=|P|$, and $b_P(u)=o_P(u)$ for
$0\le u\le\ell$. The parent base is $b=(\ell-2)/2$, the child
base is $a_k=\ell-1$, and $s_k-a_k=s'$. The child is
$\Sm(P)$, of length $2\ell-2$, and the range to be checked is
\[
m=s_k+t,\qquad 1\le t\le2\ell-2.
\]
The parent begins in zero, so the ones in a prefix of length $j$ of
its interior number $b_P(j+1)$. Counting the ones read from the two
tapes using \eqref{eq:read-counters} gives
\begin{equation}\label{eq:odd-prefix}
o(t)=b_P(z_-+1)+b_P(o_-+1).
\end{equation}
On the parent even block,
\[
G(s'+u)=b+b_P(u+1),\qquad F(s'+u)=b+u-b_P(u)
\quad(0\le u\le\ell-1).
\]
The first row of the table and \eqref{eq:odd-prefix} show that the
right side of \eqref{eq:Frec} equals
$2b+o(t)+1=a_k+o(t)=F(m)$.

For \eqref{eq:Grec}, the required identity is
\begin{equation}\label{eq:odd-park}
b_P(z)+b_P(o_-+1)=o_-.
\end{equation}
If the $t$-th output is zero, then $z=z_-+1$ and $o(t)=o_-$,
so this is exactly \eqref{eq:odd-prefix}. If the output is one,
then $z=z_-$ and $o(t)=o_-+1$. The parked-head observation says
that $P[z]=1$, with zero-based indexing; hence
$b_P(z+1)-b_P(z)=1$. Subtracting this one from
\eqref{eq:odd-prefix} proves \eqref{eq:odd-park}. The parent
$F$-formula now gives
\[
F(s'+z)+F(s'+o_-+1)
=2b+z+o_-+1-o_-=a_k+z=G(m).
\]

All these uses of the parent formulas lie in their stated domains.
Since the child has $\ell-1$ zeros and ends in one,
$0\le o_-\le\ell-2$ and $0\le z_-,z\le\ell-1$.
Thus the four parent coordinates belong to $[0,\ell-1]$.
At $t=1$, the two argument pairs are $(0,0)$ and $(1,1)$;
at $t=2\ell-2$, they are $(\ell-2,\ell-1)$ and
$(\ell-1,\ell-1)$. In particular the right endpoint of the
parent block is included, and no value beyond it is used.

\smallskip\noindent\emph{Even $k\ge2$.}
The parent $P=U_{k-1}$ now has length $\ell$, base $b=\ell/2$,
and block length $\ell$. The child base is $a_k=\ell$, and
$s=s_k=s'+\ell$. Its word $\Sp(P)$ has length $2\ell+2$, but
its core block has length $2\ell+1$. Consequently we must check
\[
m=s+t,\qquad 1\le t\le2\ell+1.
\]
Define an auxiliary prefix count only at the additional tape position:
\[
\widehat b_P(u)=b_P(u)\quad(0\le u\le\ell),\qquad
\widehat b_P(\ell+1)=b_P(\ell)+1.
\]
This is the count for the parent with its appended terminal one; it
does not extend the sequence to any undefined index. For every
$1\le v\le2\ell+2$, counting ones in the prefixes of $0P1$ and
$P$ using \eqref{eq:read-counters} gives
\begin{equation}\label{eq:even-prefix}
o(v)=\widehat b_P(z(v-1))+b_P(o(v-1)).
\end{equation}
The parent formulas are
\[
F(s'+u)=b+b_P(u),\qquad G(s'+u)=b+u-b_P(u)
\quad(0\le u\le\ell).
\]

For the $F$-identity we use
\begin{equation}\label{eq:even-park}
b_P(z-1)+b_P(o_-)=o_-.
\end{equation}
If the $t$-th output is zero, this follows from
\eqref{eq:even-prefix} because $z-1=z_-$ and $o(t)=o_-$.
If it is one, the last consumed symbol of $0P1$, at position $z$,
is one. For $z\le\ell$ this says $b_P(z)-b_P(z-1)=1$; for
$z=\ell+1$ it is the appended one and says
$\widehat b_P(\ell+1)-b_P(\ell)=1$. Subtracting that one from
\eqref{eq:even-prefix} again gives \eqref{eq:even-park}.
It follows that
\[
G(s'+z-1)+G(s'+o_-)+1
=\ell+z+o_--b_P(z-1)-b_P(o_-)=\ell+z=F(m).
\]
Here $1\le z\le\ell+1$ and $0\le o_-\le\ell$, so both
parent coordinates lie in $[0,\ell]$.

For the $G$-identity, $0\le o(t)\le\ell$ because $t$ is a proper
prefix of a word ending in one. The other coordinate $z(t)$ may
equal $\ell+1$. At that single additional parent coordinate, the
word definition in the child block, whose first bit is zero, gives
\begin{equation}\label{eq:extra-parent-point}
F(s'+\ell+1)=F(s+1)=\ell+1=b+\widehat b_P(\ell+1).
\end{equation}
Thus $F(s'+u)=b+\widehat b_P(u)$ holds also for $u=\ell+1$.
There is no appeal to the recursion being proved: both sides of
\eqref{eq:extra-parent-point} are already fixed by the word
construction. Moreover, the extra point is used only if
$z(t)=\ell+1$, which implies $t\ge\ell+1>1$; it is therefore
strictly earlier than the current index $s+t$.
Equation~\eqref{eq:even-prefix} at $v=t+1$ now gives
\[
F(s'+o(t))+F(s'+z(t))
=\ell+b_P(o(t))+\widehat b_P(z(t))
=\ell+o(t+1)=G(m).
\]
At $t=1$, the argument pairs are $(0,0)$ and $(0,1)$.
At $t=2\ell+1$, they are $(\ell,\ell-1)$ and
$(\ell,\ell+1)$: the terminal run of $\Sp(P)$ contains at least
two ones. At this final point, \eqref{eq:even-prefix} is used at
the full child length $v=2\ell+2$, where it remains valid. The
argument also covers $k=2$: its parent $(10)^{p-1}01$ need not begin in zero
for any of the even-case identities. This completes the verification
at every index.
\end{proof}

\subsection{Proof of the global family theorem}
\begin{proof}[Proof of Theorem~\ref{thm:family-global}]
The balanced initial word ends in one, so its first plus transform is
causal by Lemma~\ref{lem:maps}. Subsequent words begin in zero and end
in one; both maps remain causal. Direct tape tracing, valid for any
number of copies of $10$, gives
\[
 U_{p,2}=(0110)^{p-1}001011.
\]
Repeated use of Lemma~\ref{lem:factor-transport} proves
\eqref{eq:family-factor}. The complemented factors advance by the
opposite map, exactly as required by the parity of their level.
Primitivity of the elementary factors gives the excursion count.
Their lengths give $L_{p,k}$ in \eqref{eq:family-skeleton}.

The core construction and \eqref{eq:core-bound} establish positivity,
oddness, the bound by $n$, and the two-step law after reconstruction.
Each core is unbounded, nondecreasing, and has unit steps, so it assumes
every nonnegative value finitely often. At odd and even sequence
indices, substitution of \eqref{eq:family-lift} in the nested recursion
is respectively \eqref{eq:Frec} and \eqref{eq:Grec}.
Proposition~\ref{thm:cores} therefore proves the recurrence for $n>6p+4$.
The reconstructed values satisfy $1\le Q_p(n)\le n$, which makes every
recursive argument positive and strictly backward.

It remains to justify the shorter seeds. Formula~\eqref{eq:family-seed}
follows from the prefix of $(10)^{p-1}01$. For $p\ge3$ the following
eleven substitutions verify all missing indices $n=6p+b$; both
arguments are listed relative to $2p$ and the result relative to $4p$:
\begin{center}
\begin{tabular}{@{}rrrr@{}}\toprule
$b$ & first argument $-2p$ & second argument $-2p$ & $Q_p(n)-4p$\\\midrule
$-6$&$-3$&$-1$&$-5$\\
$-5$&$0$&$-2$&$-3$\\
$-4$&$-1$&$1$&$-3$\\
$-3$&$0$&$0$&$-1$\\
$-2$&$-1$&$1$&$-3$\\
$-1$&$2$&$0$&$-1$\\
$0$&$1$&$3$&$-1$\\
$1$&$2$&$2$&$-1$\\
$2$&$3$&$3$&$1$\\
$3$&$2$&$4$&$1$\\
$4$&$3$&$3$&$1$\\\bottomrule
\end{tabular}
\end{center}
All entries follow from the two explicitly known initial words. The
needed small argument values are $Q_p(2p+a)=2p-3$ for $a=-3,-2$,
$2p-1$ for $a=-1,0,1,2$, and $2p+1$ for $a=3,4$.
Every argument is positive and smaller than $n$. For $p=2$, the remaining initial identities are the following direct
substitutions in the independently defined cores:
\begin{center}
\begin{tabular}{@{}ccccc@{}}\toprule
$m$ & $F_2$-arguments & $G_2$-arguments & $F_2(m)$ & $G_2(m)$\\\midrule
4 & --- & $(2,3)$ & 2 & 2\\
5 & $(2,2)$ & $(2,3)$ & 3 & 2\\
6 & $(3,2)$ & $(3,4)$ & 3 & 3\\
7 & $(3,3)$ & $(4,4)$ & 3 & 4\\
8 & $(3,4)$ & $(4,4)$ & 4 & 4\\\bottomrule
\end{tabular}
\end{center}
They verify all sequence indices $8\le n\le16$; the general block proof
handles every later one.
The independently constructed sequence satisfies the prescribed
recurrence, and strong induction proves uniqueness.
\end{proof}

\begin{proposition}[Exact symmetric mode]\label{prop:sum}
For $m=s_{p,k}+t$ in the corresponding core block,
\begin{equation}\label{eq:symmetric-mode}
 F_p(m)+G_p(m)-m
 =d_{p,k}+\mathbf1_{\{k\text{ even}\}}U_{p,k}[t].
\end{equation}
Thus $Q_p(2m-1)+Q_p(2m)=2m+O_p(\log m)$ and
\begin{equation}\label{eq:pair-cancellation}
 \left(\frac{Q_p(2m-1)}{2m-1}-\frac12\right)
 +\left(\frac{Q_p(2m)}{2m}-\frac12\right)
 =O_p\left(\frac{\log m}{m}\right).
\end{equation}
\end{proposition}
\begin{proof}
Use $o(t)+z(t)=t$ and $o(t+1)=o(t)+U_{p,k}[t]$ in
\eqref{eq:family-core-odd}--\eqref{eq:family-core-even}.
Since $s_{p,k}\asymp_p2^k$, the offset is $O_p(\log m)$.
Replacing the odd denominator by $2m$ costs $O(1/m)$ because
$Q_p(2m-1)\le2m-1$.
\end{proof}

\section{The elementary Gaussian profile}\label{sec:gaussian}

Let $\varphi$ and $\Phi$ be the standard normal density and distribution
function, and define
\begin{equation}\label{eq:J}
J(x)=\varphi(\Phi^{-1}(x))\quad(0<x<1),\qquad J(0)=J(1)=0.
\end{equation}
It is a quantile-coordinate profile, not a Gaussian density in $x$.
Elementary differentiation gives
\begin{equation}\label{eq:J-ode}
J'(x)=-\Phi^{-1}(x),\qquad J(x)J''(x)=-1.
\end{equation}
The function is symmetric, strictly positive inside the interval,
strictly concave, and maximized at $1/2$.

\begin{theorem}[Uniform elementary-arch limit]\label{thm:gaussian}
With linear interpolation of integer-prefix heights,
\begin{equation}\label{eq:gaussian}
\sup_{0\le x\le1}
\left|\frac{\sqrt{k}}{\ell_k}\,h_{A_k}(\ell_k x)-J(x)\right|
\longrightarrow0.
\end{equation}
In particular,
\begin{equation}\label{eq:height-max}
\max_t h_{A_k}(t)\sim\frac{\ell_k}{\sqrt{2\pi k}}.
\end{equation}
All maximizing positions, divided by $\ell_k$, tend to $1/2$.
\end{theorem}

The $\pm2$ changes in word length must not be dropped silently. A small
exact coordinate adjustment removes them and the additive height term
simultaneously.

\subsection{An exact shifted-height identity}

Put
\begin{equation}\label{eq:shifted}
B_k=\frac{2^k}{3},\qquad \eps_k=\frac{(-1)^k}{3},\qquad
\ell_k=B_k+2\eps_k,\qquad
v_k(x)=\frac{h_{A_k}(B_kx+\eps_k)-\eps_k}{B_k}.
\end{equation}
When needed, extend a height linearly by $h(t)=t$ just to the left of
zero and by $h(t)=\ell_k-t$ just to the right of $\ell_k$.
Then $v_k$ is defined on $[0,1]$, is nonnegative and $1$-Lipschitz,
and satisfies $v_k(0)=v_k(1)=0$.
We also retain its actual word grid
\[
\mathcal G_k=\{(t-\eps_k)/B_k:0\le t\le\ell_k,\ t\in\mathbb Z\},
\]
which extends a distance $2^{-k}$ beyond $[0,1]$ when $k$ is even.
At such external grid points $v_k$ is negative of magnitude $2^{-k}$.

\begin{lemma}[Exact two-head identity]\label{lem:averaging}
At a nonempty output prefix of the map from $A_k$ to $A_{k+1}$, let
$x\in\mathcal G_{k+1}$ be its coordinate, $y=v_{k+1}(x)$, and $c$ its
last bit. Then
\begin{equation}\label{eq:averaging}
y=\frac{v_k(x+y+\delta)+v_k(x-y-\delta)}2,
\qquad
\delta=\frac{2c-1-\eps_k}{2B_k},\quad |\delta|\le2^{1-k}.
\end{equation}
The parent arguments are shifted word-grid points, with the single
possible extra terminal point inherited from $\Sp$.
If a parent argument is outside $[0,1]$, the output satisfies $y=1-x$.
The zero-length and complete-prefix boundary cases are linear pieces
with $y=x$ or $y=1-x$ when their coordinates lie in $[0,1]$.
\end{lemma}
\begin{proof}
Write $H$ for the raw output height and $t$ for its prefix length.
For $\Sp$ let $r=i-1$, $u=j$; for $\Sm$ let $r=i$, $u=j+1$.
With $\sigma=(-1)^{k+1}$, the tape heights give
\[
H=h_{A_k}(r)+h_{A_k}(u)+\sigma,\quad
r+u=t-\sigma,\quad r-u=H+2c-1.
\]
For the extra terminal prefix in $\Sp$, set $h_{A_k}(\ell_k+1)=-1$.
Since $\sigma=-3\eps_k$, $\eps_{k+1}=-\eps_k$, and $B_{k+1}=2B_k$,
substitution in \eqref{eq:shifted} gives \eqref{eq:averaging} exactly.

No left external parent occurs at a nonempty prefix: for $\Sm$,
$r,u\ge1$, whereas for $\Sp$, $r,u\ge0$ and $\eps_k<0$.
Writing $B=B_k$ and $\ell=\ell_k$, the complete list of right external
parents is
\begin{equation}\label{eq:external-parents}
\begin{array}{c|c|c|c}
 \text{map}&\text{raw parent index }r& (r-\eps_k)/B
   &v_k((r-\eps_k)/B)\\ \hline
 \Sm\ (k\text{ even})&\ell&1+1/(3B)&-1/(3B)\\
 \Sp\ (k\text{ odd})&\ell+1&1+2/(3B)&-2/(3B).
\end{array}
\end{equation}
The other parent $u$ lies in the internal shifted grid.
In either row the head identities give
$\ell-u=h_{A_k}(u)+2c-2$.
Since $h_{A_k}(u)\le\ell-u$, they force $c=1$ and
$h_{A_k}(u)=\ell-u$. Thus all remaining symbols on the other tape are
ones, and $H=\ell_{k+1}-t$, or $y=1-x$.
The empty and complete output prefixes belong to $[0,1]$ exactly when
$k$ is even; their $(x,y)$ coordinates are respectively
\begin{equation}\label{eq:internal-prefix-endpoints}
 \left(\frac1{6B},\frac1{6B}\right),\qquad
 \left(1-\frac1{6B},\frac1{6B}\right).
\end{equation}
For odd $k$ the corresponding coordinates are outside $[0,1]$.
\end{proof}

\subsection{Upper and positive lower barriers}

\begin{lemma}[Coarse upper bound]\label{lem:upper}
There is an absolute constant $C$ such that
\begin{equation}\label{eq:upper}
0\le v_k(x)\le C/\sqrt{k}\qquad(0\le x\le1).
\end{equation}
\end{lemma}
\begin{proof}
Use the quadratic $q(x)=1+x(1-x)$ and amplitudes
$b_4=1$, $b_{k+1}=b_k/(1+b_k^2/32)$.
We prove $v_k\le b_kq$ on the shifted grids. The initial bound holds
because $v_4\le1/2$ on $[0,1]$, while the external grid values are
negative. Also
\[
b_{k+1}^{-2}-b_k^{-2}=\frac1{16}+\frac{b_k^2}{1024},
\]
so $b_k\asymp k^{-1/2}$, and $2^{1-k}\le b_k/4$ for $k\ge4$.

Suppose the next bound fails at an internal output grid point. Write
$b=b_k$, $b'=b_{k+1}$. The grid lies within $2^{-k-1}$ of $[0,1]$, so
$q\ge0.9$ there. Then
$y>b'q(x)>0.87b$, hence $|y+\delta|>b/2$.
The exact quadratic average and \eqref{eq:averaging} give
\[
y\le b\{q(x)-(y+\delta)^2\}<bq(x)-b^3/4<b'q(x),
\]
since $q(x)\le5/4$ and $b-b'\le b^3/32$.
The extra terminal parent value is negative and satisfies the same
upper comparison; $q$ is positive on all points in question.
The output boundary points have height $O(2^{-k})$ and satisfy the bound
directly. Interpolation adds at most $O(2^{-k})$, proving \eqref{eq:upper}.
\end{proof}

\begin{lemma}[Positive sine barrier]\label{lem:sine}
Let $b_6=1/64$ and
\begin{equation}\label{eq:sine-rec}
b_{k+1}=b_k\cos\!\bigl(\pi(b_k+2^{1-k})\bigr),\qquad k\ge6.
\end{equation}
At every point $x\in\mathcal G_k\cap[0,1]$,
\begin{equation}\label{eq:sine}
v_k(x)\ge b_k\sin(\pi x),\qquad
\sqrt{k}\,b_k\longrightarrow1/\pi.
\end{equation}
\end{lemma}
\begin{proof}
At level six, each internal raw height is at least one; after the shift
its normalized height is at least $1/32$. Thus the initial grid bound
holds. The amplitudes stay positive and decrease, with
$b_k+2^{1-k}<1/8$.
A failure $0\le y<b_{k+1}\sin(\pi x)$ at the next grid point, when both
parents lie in $[0,1]$, contradicts the sine addition formula:
\[
y\ge b_k\sin(\pi x)\cos\bigl(\pi(y+\delta)\bigr)
\ge b_k\sin(\pi x)\cos\bigl(\pi(b_k+2^{1-k})\bigr).
\]
If a parent lies outside, Lemma~\ref{lem:averaging} gives
$y=1-x\ge\sin(\pi x)/\pi$, excluding a failure. The linear boundary
cases are treated by the same inequality and its reflection.

For the asymptotic let $c_k=b_k^{-2}$.
The amplitudes tend to zero, since a positive limit in
\eqref{eq:sine-rec} would give a fixed strict contraction. Moreover
\[
c_{k+1}-c_k=c_k\tan^2\!\bigl(\pi(b_k+2^{1-k})\bigr)
\le8\pi^2(1+c_k4^{1-k}).
\]
Iteration, using convergence of the corresponding infinite product,
gives $c_k=O(k)$. Hence $2^{1-k}/b_k\to0$ and the exact difference tends
to $\pi^2$. Averaging yields $c_k/k\to\pi^2$.
\end{proof}

\subsection{Gaussian comparison and uniform convergence}

\begin{lemma}[Smooth comparison functions]\label{lem:smooth}
There are smooth functions of the following two types arbitrarily close
to $J$, uniformly on $[0,1]$.
For the upper type, $f>0$ on a neighborhood of $[0,1]$ and $ff''=-1$.
For the lower type, $f$ is concave, negative at both endpoints, positive
on a central interval compactly contained in $(0,1)$, and
$ff''\ge-1$ wherever $f>0$; it is affine near its two zeros.
For the lower type it is the positive part $f_+$ which approximates $J$.
\end{lemma}
\begin{proof}
Upper functions are
$f_\eta(x)=(1+2\eta)J((x+\eta)/(1+2\eta))$, $\eta>0$.
For lower functions fix $0<A<1/\sqrt{2\pi}$ and a smooth $g_\eta$ which
is zero for $v\le\eta$, equals $1/v$ for $v\ge2\eta$, and satisfies
$0\le g_\eta(v)\le1/v$ when $v>0$.
Solve symmetrically
\[
f''=-g_\eta(f),\qquad f(1/2)=A,\qquad f'(1/2)=0,
\]
and continue affinely once $f\le\eta$.
Above height $2\eta$ this agrees with the solution of $ff''=-1$ of
central height $A$. That solution has half-width $A\sqrt{\pi/2}$,
since $|f'|=\sqrt{2\log(A/f)}$ and
$\int_0^A(2\log(A/v))^{-1/2}\,dv=A\sqrt{\pi/2}$.
The remaining time from height $2\eta$ to zero is at most
$2\eta/\sqrt{2\log(A/(2\eta))}$.
Thus the zeros lie inside $(0,1)$ for small $\eta$.
The positive part differs by at most $2\eta$ from the compactly supported
Gaussian-quantile arch of height $A$. Let $A$ tend upward to
$1/\sqrt{2\pi}$ and then $\eta$ tend to zero.
\end{proof}

\begin{lemma}[One-step comparison]\label{lem:comparison}
Fix an upper-type function $f$ and $\gamma>1$, or a lower-type function
$f$ and $0<\gamma<1$. There is $D_0=D_0(f,\gamma)\ge4$ such that,
whenever $D\ge D_0$ and $B_k\ge D^3$, the corresponding implication
holds on the internal shifted grids:
\begin{align}
 v_k\le\frac{\gamma f}{\sqrt D}
 &\quad\Longrightarrow\quad
 v_{k+1}\le\frac{\gamma f}{\sqrt{D+1}}
 &&\text{(upper type)},\label{eq:comparison-upper}\\
 v_k\ge\frac{\gamma f}{\sqrt D}
 &\quad\Longrightarrow\quad
 v_{k+1}\ge\frac{\gamma f}{\sqrt{D+1}}
 &&\text{(lower type)}.\label{eq:comparison-lower}
\end{align}
Here the parent inequality is assumed on $\mathcal G_k\cap[0,1]$,
and the conclusion is on $\mathcal G_{k+1}\cap[0,1]$.
\end{lemma}
\begin{proof}
Put $B=B_k$, $a=\gamma D^{-1/2}$, $b=\gamma(D+1)^{-1/2}$, and
$g=|\gamma^2-1|>0$. We use the explicit scale comparison
\begin{equation}\label{eq:comparison-strict-margin}
 \frac ba-\left(1-\frac{b^2}{2}\right)
 =\frac{\gamma^2-1}{2D}+\varepsilon_D,
 \qquad |\varepsilon_D|\le\frac{C_\gamma}{D^2},
 \qquad C_\gamma=\frac{g+1}{2},\quad D\ge4.
\end{equation}
Indeed, with $s=(D+1)^{-1}$, Taylor's theorem gives
\[
 \left|\sqrt{1-s}-1+\frac s2\right|
 \le\frac{s^2}{8(1-s)^{3/2}}\le\frac{s^2}{4},
\]
and replacing $(D+1)^{-1}$ by $D^{-1}$ gives the displayed estimate.

For a nonempty output prefix write $y=v_{k+1}(x)$ and
$d=y+\delta$, where $|\delta|\le B^{-1}$ by
Lemma~\ref{lem:averaging}. Thus $y$ is the mean of the two parent
values at $x\pm d$. The constants below depend only on $f,\gamma$,
and are uniform in $D,B,x$, and the choice of word map.

For clarity, the Taylor estimate used in both directions is
\begin{equation}\label{eq:symmetric-taylor-bound}
 \frac{f(x+d)+f(x-d)}2
 =f(x)+\frac12f''(x)d^2+R(x,d),\qquad
 |R(x,d)|\le\frac{d^4}{24}
    \max_{|z-x|\le|d|}|f^{(4)}(z)|.
\end{equation}
It follows by applying fourth-order remainder bounds separately at
$x+d$ and $x-d$; the odd terms cancel. The whole sampling segment lies
in $I$ in the upper case and in $[0,1]$ in the lower case considered below.

For an upper comparison choose a compact interval $I$ in the positive
neighborhood on which $ff''=-1$, with $[0,1]$ in its interior, and put
\[
 m=\min_I f>0,\qquad M=\max_I f,\qquad
 B_4=\max_I|f^{(4)}|.
\]
All possible parent arguments lie within $B^{-1}$ of $[0,1]$.
Increase a preliminary threshold $D_{\rm geom}\ge4$ so that these
arguments lie in $I$ for $D\ge D_{\rm geom}$ and $B\ge D^3$.
Every external parent value is negative, including the additional
terminal value of $\Sp$, so it automatically satisfies the upper
comparison. The parent inequality therefore gives $0\le y\le aM$.
The symmetric Taylor expansion has remainder at most $B_4d^4/24$.
At a putative failure $y>bf(x)$ we have
\[
 d^2\ge b^2f(x)^2-\frac{2aM}{B},\qquad
 |d|\le\frac{\gamma M+1}{\sqrt D}.
\]
Since $f''=-1/f$, it follows that
\begin{align*}
 y&\le af(x)\left(1-\frac{b^2}{2}+E_+\right),\\
 0\le E_+&\le\frac{\gamma M}{m^2\sqrt D\,B}
       +\frac{B_4(\gamma M+1)^4}{24mD^2}
       \le\frac{C_+}{D^2},\\
 C_+&=\frac{\gamma M}{m^2}
       +\frac{B_4(\gamma M+1)^4}{24m}.
\end{align*}
Choose $D_0\ge D_{\rm geom}$ also so that
$D_0\ge4(C_++C_\gamma)/g$. Equation~\eqref{eq:comparison-strict-margin}
then gives
\[
 1-\frac{b^2}{2}+E_+\le\frac ba-\frac{g}{4D},
\]
contradicting $y>bf(x)$.

For a lower comparison a failure can occur only where $f(x)>0$,
since $v_{k+1}\ge0$ on $[0,1]$. Let $(p,q)\Subset(0,1)$ be the
positive set of $f$, and set $M=\max_{[0,1]}f_+$.
If a parent lies outside $[0,1]$, Lemma~\ref{lem:averaging} gives
$y=1-x\ge1-q$, exceeding $bM$ after increasing $D_{\rm geom}$.
Thus a lower failure can involve only internal parent points.

Choose disjoint affine neighborhoods of the two zeros of $f$, and
then choose $\rho>0$ so that the compact set
\[
 K_\rho=\{x\in[0,1]:0\le f(x)\le\rho\}
\]
lies in the interior of their union. Let $\eta>0$ be its distance
from the complement of that union. Increase $D_{\rm geom}$ until
$b\rho+B^{-1}<\eta$. If a putative failure has $0<f(x)\le\rho$,
then $|d|\le b\rho+B^{-1}<\eta$, so the entire sampling segment
lies in the same affine neighborhood. Its exact average is $f(x)$,
giving $y\ge af(x)>bf(x)$, a contradiction. This includes sampling
across a zero and requires no division by a vanishing function.

It remains to consider $f(x)\ge\rho$. Set
$B_j=\max_{[0,1]}|f^{(j)}|$ for $j=2,4$. Now
\[
 d^2\le b^2f(x)^2+\frac{2bM}{B}+B^{-2},\qquad
 |d|\le\frac{\gamma M+1}{\sqrt D}.
\]
Using $f''\le0$ and $ff''\ge-1$, symmetric Taylor expansion gives
\begin{align*}
 y&\ge af(x)\left(1-\frac{b^2}{2}-E_-\right),\\
 0\le E_-&\le\frac{B_2}{\rho}
       \left(\frac{\gamma M}{\sqrt D\,B}+\frac1{2B^2}\right)
       +\frac{B_4(\gamma M+1)^4}{24\rho D^2}
       \le\frac{C_-}{D^2},\\
 C_-&=\frac{B_2(\gamma M+1/2)}{\rho}
       +\frac{B_4(\gamma M+1)^4}{24\rho}.
\end{align*}
Choose $D_0\ge D_{\rm geom}$ also so that
$D_0\ge4(C_-+C_\gamma)/g$. Since $\gamma^2-1=-g$,
Equation~\eqref{eq:comparison-strict-margin} gives
\[
 1-\frac{b^2}{2}-E_-\ge\frac ba+\frac{g}{4D},
\]
again contradicting the proposed failure.

Finally, the only internal empty or complete prefixes are the two
points in \eqref{eq:internal-prefix-endpoints}. For the upper comparison
their height satisfies
\[
 \frac1{6B}\le\frac1{6D^3}
 \le\frac{\gamma m}{\sqrt{D+1}}\le bf(x)
\]
after a fixed increase of $D_{\rm geom}$. For the lower comparison
choose $\varepsilon_0>0$ such that $f<0$ on
$[0,\varepsilon_0]\cup[1-\varepsilon_0,1]$; it suffices that
$(6D^3)^{-1}\le\varepsilon_0$.
All geometric choices above are uniform: if
$\tau=\operatorname{dist}([0,1],\mathbb R\setminus I)>0$, it is enough
in the upper case to require $D^{-3}<\tau$ and the displayed endpoint
inequality. In the lower case the other requirements follow from
\[
 \frac{\gamma M}{\sqrt D}<1-q,\qquad
 \frac{\gamma\rho}{\sqrt D}+D^{-3}<\eta.
\]
These are finitely many conditions involving only the fixed function
$f$ and multiplier $\gamma$. Together with
$D_0\ge4(C_\pm+C_\gamma)/g$, they specify a threshold independent of
$k$, $B$, and the output point. This completes both implications on
the stated internal grids.
\end{proof}

\begin{proof}[Proof of Theorem~\ref{thm:gaussian}]
Fix upper and lower functions $f^+,f^-$ and parameters
$\gamma_+>1$ and $0<\gamma_-<1$. We first initialize the comparisons
at one sufficiently large level $K$.

By Lemma~\ref{lem:upper}, $v_K\le C_0/\sqrt K$ on $[0,1]$.
Choose a fixed $c>0$ so small that
$\gamma_+\min_{[0,1]}f^+/\sqrt c\ge C_0$.
Then the upper comparison holds with $D_+=cK$.
For the lower comparison, the positive support of $f^-$ is a compact
subinterval of $(0,1)$, so
\[
 R=\max_{x:\,f^-(x)>0}\frac{f^-(x)}{\sin(\pi x)}<\infty.
\]
Lemma~\ref{lem:sine} gives $v_K\ge c_0\sin(\pi x)/\sqrt K$ on the
internal grid for some fixed $c_0>0$ and all large $K$.
Choose a fixed $C>0$ with $\gamma_-R/\sqrt C\le c_0$.
The lower comparison then holds with $D_-=CK$; where $f^-\le0$
it follows from $v_K\ge0$.

Choose $K$ large enough that both $D_\pm\ge D_0$ and
$B_K\ge\max(D_+^3,D_-^3)$. These inequalities persist under iteration:
$B_k$ doubles while $D$ increases by one, and
$2D^3\ge(D+1)^3$ for $D\ge4$. Thus Lemma~\ref{lem:comparison}
propagates the two bounds, with time parameters $D_\pm+k-K$.
Since
\[
 \frac{\sqrt k}{\sqrt{D_\pm+k-K}}\longrightarrow1,
\]
the limiting lower and upper bounds are
$\gamma_-(f^-)_+$ and $\gamma_+f^+$ uniformly on the internal grids.
Let the approximants tend uniformly to $J$ and the two multipliers
tend to one. The grid spacing is $B_k^{-1}$, and the omitted endpoint
intervals have length at most $B_k^{-1}$. The $1$-Lipschitz property
and $v_k(0)=v_k(1)=0$ therefore extend the convergence uniformly to
$[0,1]$.

Finally, $B_k/\ell_k=1+O(2^{-k})$. In normalized coordinates the
horizontal and vertical shifts in \eqref{eq:shifted} are both
$O(2^{-k})$; the height function is $1$-Lipschitz. These changes
vanish after multiplication by $\sqrt k$, and yield
\eqref{eq:gaussian}. Uniform convergence and the unique strict maximum
of $J$ imply both \eqref{eq:height-max} and convergence of every
maximizing position to $1/2$.
\end{proof}

\begin{corollary}[Height moments]\label{cor:moments}
For each fixed real $q>0$,
\[
\frac{k^{q/2}}{\ell_k^{q+1}}
\sum_{t=0}^{\ell_k-1}h_{A_k}(t)^q
\longrightarrow\frac{(2\pi)^{-q/2}}{\sqrt{q+1}}.
\]
\end{corollary}
\begin{proof}
Uniform convergence reduces the assertion to a Riemann integral.
The substitution $x=\Phi(z)$ gives
$\int_0^1J(x)^q\,dx=\int_{\mathbb R}\varphi(z)^{q+1}\,dz$,
which has the stated value.
\end{proof}

\begin{figure}[htbp]
\centering
\includegraphics[width=.96\textwidth]{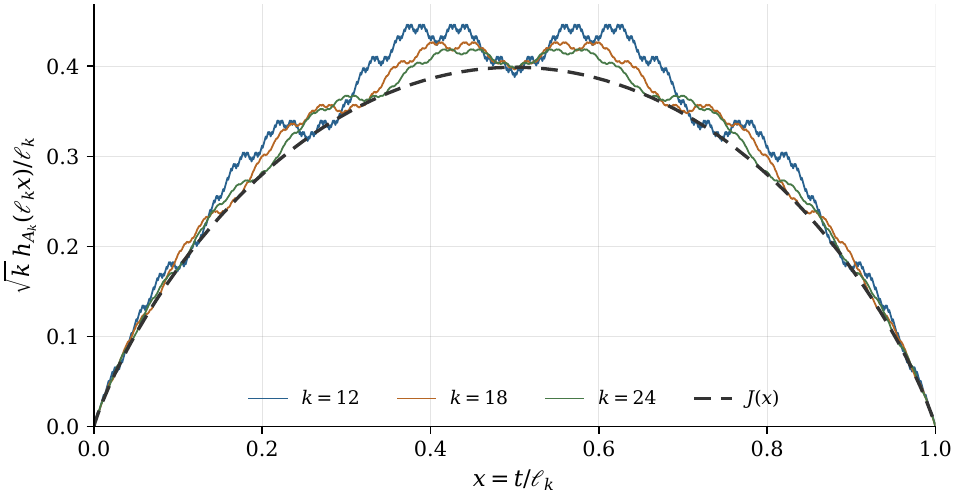}
\caption{Computed elementary-arch profiles and the proved limit.
Finite-level ripples remain visible. The normalization is specified
in Theorem~\ref{thm:gaussian}; no fitted amplitude or horizontal shift
is used.}
\label{fig:gaussian}
\end{figure}

\section{Family profiles and sharp fluctuations}\label{sec:family-profiles}
For the standard-block statements, let
\[
 N_{j,r}(t)=2\left\lfloor\frac{2^jt-r}{2}\right\rfloor+r,
 \qquad r\in\{0,1\},\quad 1\le t\le2.
\]
This samples the specified parity within two indices of $2^jt$.

The elementary words $A_k$ are the same for every $p$. Define the continuous signed height
profile $\mathscr H_p$ on $[0,c_p]$ by
\begin{equation}\label{eq:family-height-profile}
\mathscr H_p(x)=
\begin{cases}
 J(x-3i),&3i\le x\le3i+1,\quad 0\le i\le p-2,\\
 -2J((x-3i-1)/2),&3i+1\le x\le3i+3,\quad0\le i\le p-2,\\
 4J((x-3p+3)/4),&3p-3\le x\le c_p.
\end{cases}
\end{equation}
Joining formulas agree at their endpoints.

\begin{theorem}[Uniform family profiles]\label{thm:family-profile}
Fix $p\ge2$. With $B_k=2^k/3$, linear interpolation and the harmless
zero continuation over a missing endpoint interval of length $O(1)$,
\begin{equation}\label{eq:family-height-limit}
 \frac{\sqrt{k}}{B_k}h_{U_{p,k}}(B_kx)
 \longrightarrow\mathscr H_p(x)
 \quad\text{uniformly on }[0,c_p].
\end{equation}
For indices of parity $r$ whose core coordinate satisfies
$m=s_{p,k}+B_kx+O(1)$ within the $k$-th macro block,
\begin{equation}\label{eq:family-natural-ratio}
 (-1)^{k+r+1}\sqrt{k}\left(\frac{Q_p(n)}n-\frac12\right)
 \longrightarrow\frac{\mathscr H_p(x)}{2(c_p+x)}
\end{equation}
uniformly in the normalized block coordinate, with endpoints taken
from the actual core block.
On standard dyadic blocks let
\[
 \lambda_p(t)=\left\lfloor\log_2(t/C_p)\right\rfloor,
 \qquad x_p(t)=3\,2^{-\lambda_p(t)-1}t-c_p\qquad(1\le t\le2).
\]
Then the complete limiting ratio profile is
\begin{equation}\label{eq:family-standard-profile}
 \mathcal P_p(t)=
 (-1)^{\lambda_p(t)+1}
 \frac{\mathscr H_p(x_p(t))}{2(c_p+x_p(t))},
\end{equation}
and for $r=0,1$,
\begin{equation}\label{eq:family-standard-limit}
 \sup_{1\le t\le2}\left|
 (-1)^{j+r}\sqrt j\left(\frac{Q_p(N_{j,r}(t))}{N_{j,r}(t)}-\frac12\right)
 -\mathcal P_p(t)\right|\longrightarrow0.
\end{equation}
The formula at a macro boundary is independent of which side is chosen.
The profile has exactly $2p-1$ zeros in $(1,2)$, all sign transitions.
Writing $\operatorname{mant}_2(y)=y/2^{\lfloor\log_2 y\rfloor}$, their set is
\begin{equation}\label{eq:family-centers}
\begin{split}
\mathcal Z_p={}&
 \left\{\operatorname{mant}_2\!\left(\frac23(c_p+3i+s)\right):
 0\le i\le p-2,\ s=0,1\right\}\\
 &\quad\cup
 \left\{\operatorname{mant}_2\!\left(\frac23(6p-2)\right)\right\}.
\end{split}
\end{equation}
\end{theorem}
\begin{proof}
Each factor in \eqref{eq:family-factor} is an $A_{k+s}$ with
$s\in\{0,1,2\}$. Its length divided by $B_k$ tends to $2^s$, and
$\sqrt{k}/\sqrt{k+s}\to1$. Theorem~\ref{thm:gaussian}, the finite
number of factors for fixed $p$, and the vanishing of $J$ at the
endpoints prove \eqref{eq:family-height-limit}. Length and joining
corrections are $O_p(1)$ before normalization; Lipschitz heights
make their effect $O_p(\sqrt{k}2^{-k})$.

The two core formulas give, uniformly in a macro block,
\[
 \frac{Q_p(n)}n-\frac12
 =\frac{(-1)^{k+r+1}}{2m}h_{U_{p,k}}(m-s_{p,k})
  +O_p(k2^{-k}).
\]
Also $s_{p,k}=c_pB_k+O_p(k)$. This proves
\eqref{eq:family-natural-ratio}. For a standard block the macro level
is $k=j+\lambda_p(t)$ away from an $O_p(j2^{-j})$ neighborhood of
a boundary. The normalized core coordinate tends to $x_p(t)$.
The sign in \eqref{eq:family-standard-profile} is exactly the product
of the natural and standard parity compensations. Uniformity across
the small boundary intervals follows from the endpoint zeros of
$\mathscr H_p$.

The full set of limiting beginnings of elementary excursions is
$\gamma2^k$, where $\gamma$ ranges over the unnormalized constants in
\eqref{eq:family-centers}. They all lie in $[C_p,2C_p)$, are distinct,
and none is a power of two: their numerator before division by three
is not divisible by three. Thus their binary mantissas are distinct
and lie strictly between one and two. Their number is $2p-1$.
The signs alternate at internal joins; at a macro join the level
parity changes. Positivity of $J$ proves that there are no further
limiting zeros.
\end{proof}

\begin{corollary}[Moments and sharp parameter-dependent fluctuations]
\label{cor:family-envelope}
For every fixed $p\ge2$,
\[
 Q_p(n)=n/2+O_p(n/\sqrt{\log n}).
\]
Let $b_p>0$ be the unique root of
\begin{equation}\label{eq:family-root}
 b_p\left(\frac{3p-1}{2}+\Phi(-b_p)\right)=\varphi(b_p).
\end{equation}
Then
\begin{equation}\label{eq:family-envelope}
 \limsup_{n\to\infty}\sqrt{\log_2 n}
 \left|\frac{Q_p(n)}n-\frac12\right|=\Lambda_p=\frac{b_p}{2}.
\end{equation}
The signed extremes are $\pm\Lambda_p$. The largest excursion is
always the final factor $A_{k+2}$, with maximizing local coordinate
$x_* =\Phi(-b_p)$. The constants decrease strictly with $p$, and
\begin{equation}\label{eq:family-large-p}
 \Lambda_p\sim\frac1{3p\sqrt{2\pi}}\qquad(p\to\infty).
\end{equation}
The raw absolute macro-height maximum and all of its maximizing
positions $t_k$ satisfy
\[
 \max_t|h_{U_{p,k}}(t)|\sim\frac{4B_k}{\sqrt{2\pi k}},
 \qquad \frac{t_k}{L_{p,k}}\longrightarrow\frac{3p-1}{3p+1}.
\]
This maximum concerns raw word height, not the relative error;
the latter has the different local maximizer $x_*$ above.
For each fixed real $q>0$, the absolute macro-height moment is
\begin{equation}\label{eq:family-moments}
 \frac{k^{q/2}}{L_{p,k}^{q+1}}
 \sum_{t=0}^{L_{p,k}-1}|h_{U_{p,k}}(t)|^q
 \longrightarrow
 \frac{(p-1)(1+2^{q+1})+4^{q+1}}{c_p^{q+1}}
 \frac{(2\pi)^{-q/2}}{\sqrt{q+1}}.
\end{equation}
\end{corollary}
\begin{proof}
Uniform profiles imply the order bound. On a factor of width $w$,
the limiting absolute ratio error is $J(x)/(2(\alpha+x))$, where
$\alpha$ is its starting core coordinate divided by its width.
The possible parameters are
\[
 c_p+3i,\qquad\frac{c_p+3i+1}{2}\quad(0\le i\le p-2),
 \qquad\frac{6p-2}{4}.
\]
The last is strictly the smallest. For fixed $\alpha$, the derivative numerator of $J(x)/(2(\alpha+x))$
is $(\alpha+x)J'(x)-J(x)$; its derivative is
$(\alpha+x)J''(x)<0$, and it crosses zero once.
Writing $x_* =\Phi(-b_p)$ and using $J'(x_*)=b_p$ gives
\eqref{eq:family-root}; at that point the maximum is $b_p/2$.
Each index parity supplies the opposite sign. This proves
\eqref{eq:family-envelope}.
The root equation's derivative in $b$ is $(3p-1)/2+\Phi(-b)>0$.
Increasing $p$ therefore decreases its root strictly.
As $p\to\infty$, $b_p\to0$ and the root equation gives
$b_p\sim2\varphi(0)/(3p)$.
The raw maximum follows because the final width-four factor has
strictly greater limiting absolute height than any width-one or
width-two factor; its local maximizing coordinate tends to $1/2$.
For \eqref{eq:family-moments}, sum the elementary moment formula of
Corollary~\ref{cor:moments} over the factors. Their widths are
$1,2,\ldots,1,2,4$ in units of $B_k$; a width $w$ contributes
$w^{q+1}\int_0^1J(x)^q\,dx$.
\end{proof}

\begin{center}
\begin{tabular}{@{}rrrr@{}}\toprule
$p$ & number of excursions & $\Lambda_p$ & recursion starts at\\\midrule
2&3&0.067088003576468&8\\
3&5&0.044502379773681&12\\
4&7&0.033318977308329&18\\
5&9&0.026633875777645&24\\
8&15&0.016631790335322&42\\\bottomrule
\end{tabular}
\end{center}
These are numerical evaluations of \eqref{eq:family-root}, not fitted
amplitudes. $Q_3$ and $Q_4$ have the transition sets reported by
\eqref{eq:family-centers}; in increasing order they are
\[
 \mathcal Z_3=\{13/12,7/6,4/3,5/3,11/6\},
\]
\[
 \mathcal Z_4=\{13/12,7/6,4/3,17/12,19/12,5/3,11/6\}.
\]

\begin{figure}[htbp]
\centering
\includegraphics[width=.99\textwidth]{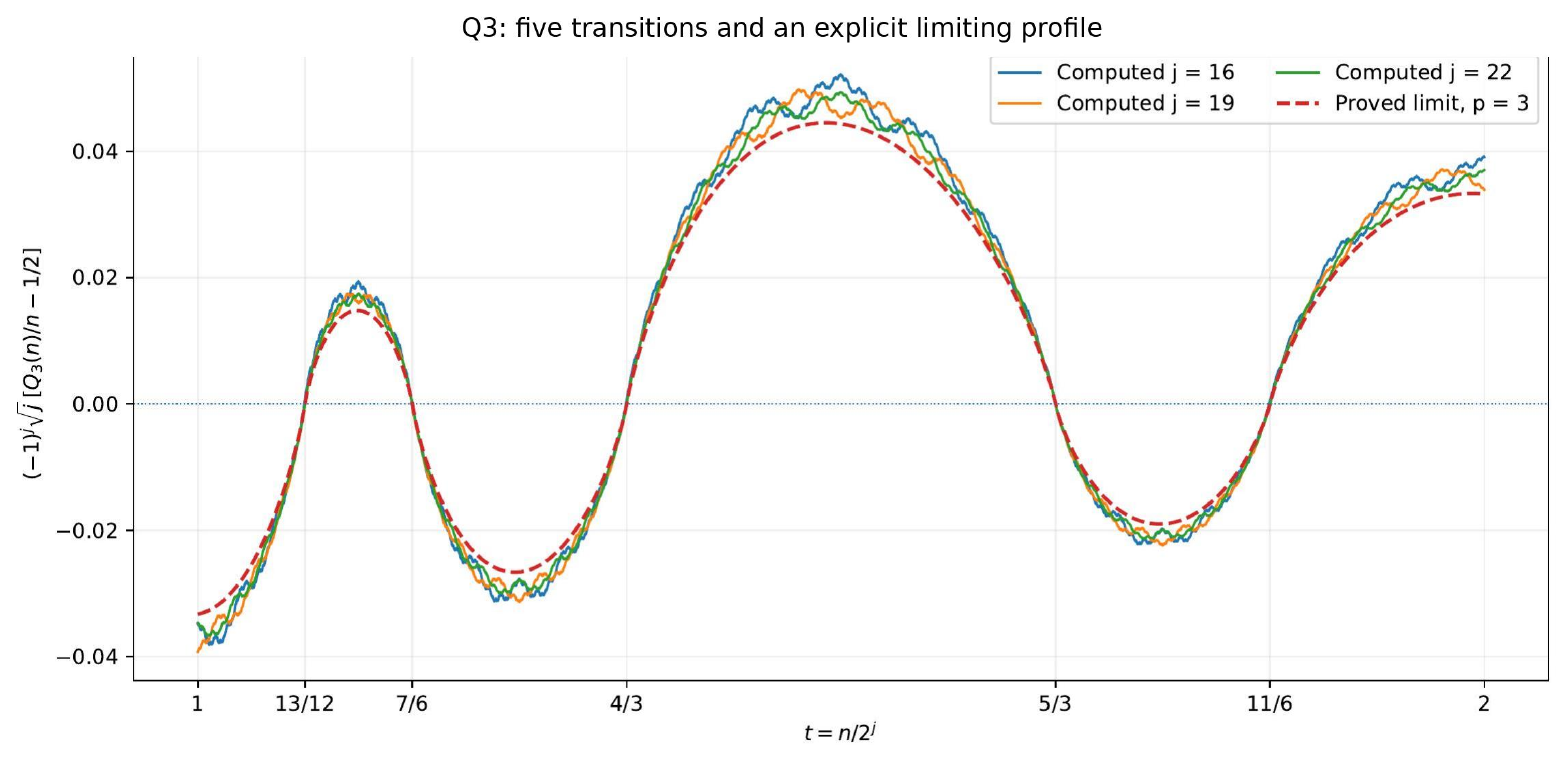}
\caption{The five-arch sequence $Q_3$: directly computed, unaveraged
even-parity profiles and the explicit $p=3$ limit
\eqref{eq:family-standard-profile}, without fitted parameters.
The sign and $\sqrt j$ factors are prescribed, not fitted.}
\label{fig:family-d5}
\end{figure}

\begin{figure}[htbp]
\centering
\includegraphics[width=.92\textwidth]{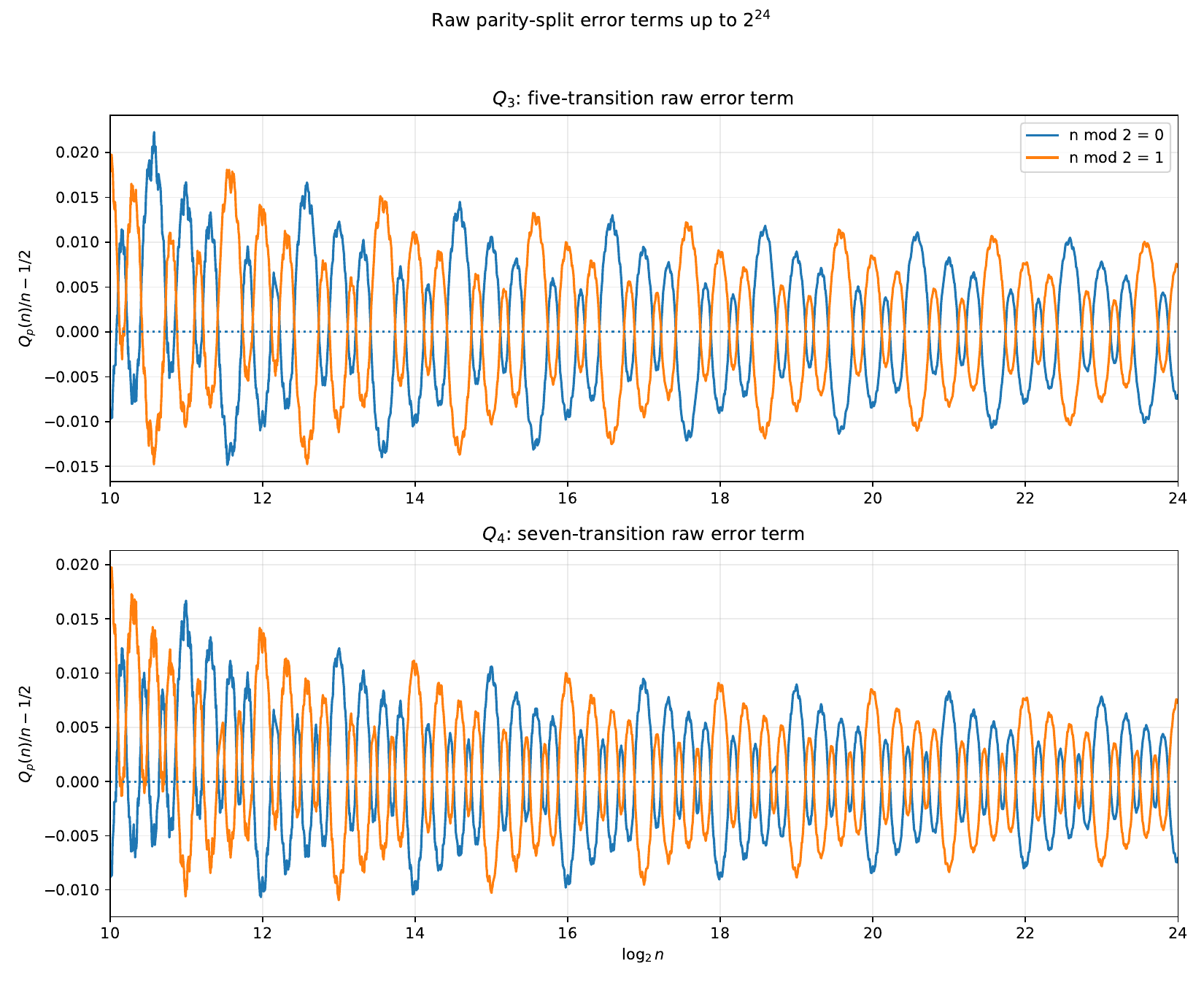}
\caption{Logarithmically sampled raw error terms $Q_p(n)/n-1/2$ for
$Q_3$ and $Q_4$, separated by parity.
These plots visualize the five- and seven-transition oscillatory
patterns directly in the unaveraged data.}
\label{fig:q34-error}
\end{figure}

Figure~\ref{fig:q34-error} records the raw parity-split error terms for
$Q_3$ and $Q_4$ on a logarithmic scale. It complements the compensated
limit plots by showing the self-similar oscillations before averaging
or sign correction.

\begin{figure}[htbp]
\centering
\includegraphics[width=.95\textwidth]{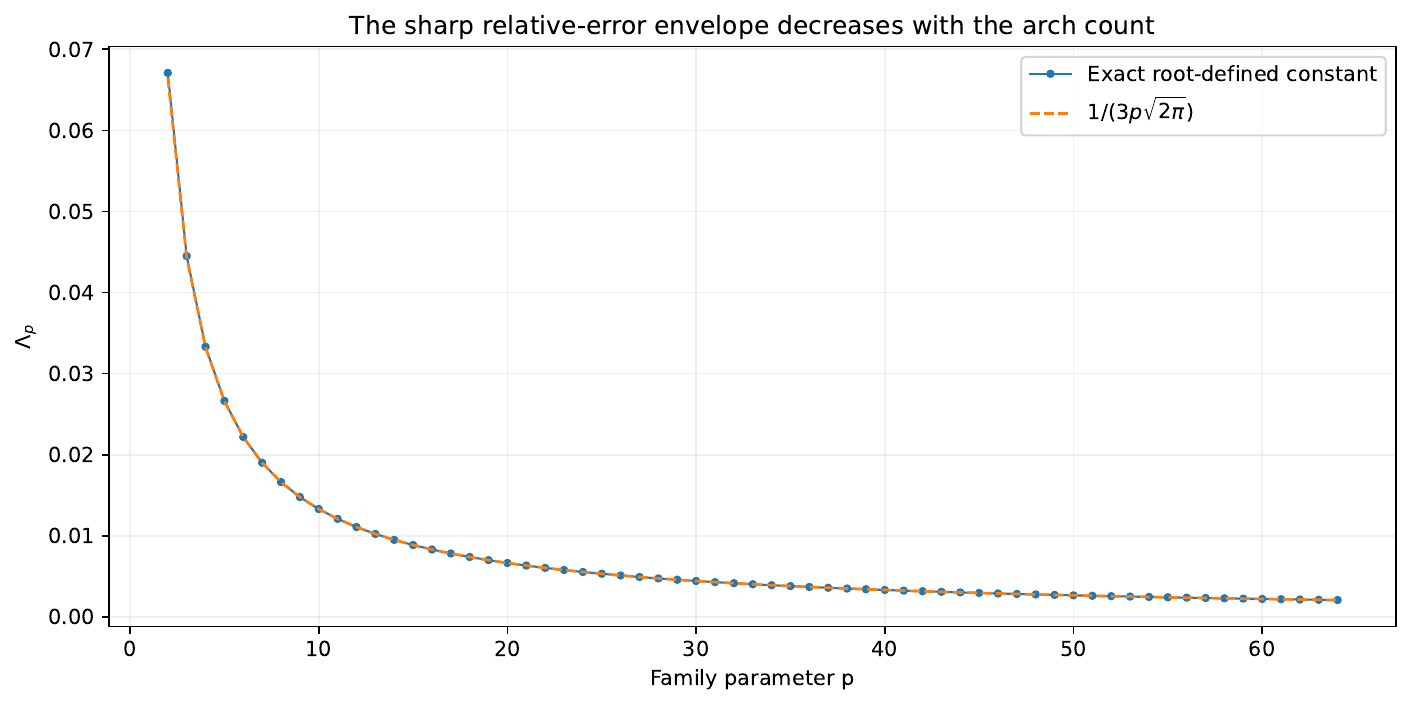}
\caption{Root-defined sharp error constants and their large-parameter
asymptotic. The approximation in $p$ is separate from the fixed-$p$
limit as the sequence index tends to infinity.}
\label{fig:family-envelope}
\end{figure}

\section{Exact finite-level transitions}\label{sec:family-boundaries}
\subsection{Low heights of the common elementary words}
The Gaussian limit controls macroscopic signs but not all finite zeros.
The following two lemmas determine the short end intervals containing
all relevant low heights, independently of any convergence rate.

\begin{lemma}[All low heights of an elementary word]\label{lem:low-heights}
Let $k\ge8$, $\lambda=\lfloor k/2\rfloor$, and $L=\ell_k$.
For integer $0\le t\le\lambda+4$,
\begin{equation}\label{eq:low-prefix}
 h_{A_k}(t)=g_\lambda(t):=
 \begin{cases}t,&0\le t\le\lambda,\\t-2,&\lambda+1\le t\le\lambda+4.
 \end{cases}
\end{equation}
The same formula holds at the other end by reflection:
$h_{A_k}(L-t)=g_\lambda(t)$. Every remaining height satisfies
\begin{equation}\label{eq:low-middle}
 h_{A_k}(t)\ge\lambda+3
 \qquad(\lambda+5\le t\le L-\lambda-5).
\end{equation}
Thus these two short end intervals contain \emph{all} positions of
height at most $\lambda+2$.
\end{lemma}
\begin{proof}
For the finite base, $A_8$ is the concatenation of the following two
strings:
\[
\begin{gathered}
\mathtt{0000100000100010110010000100010110010111011}\\
\mathtt{0010001011001011101111011001011101111101111}.
\end{gathered}
\]
Its length is $86$, its first nine heights are
$0,1,2,3,4,3,4,5,6$, and its minimum on $9\le t\le77$ is $7$.
This verifies the assertion at level eight.

Suppose it holds at level $k$. Write $\sigma=1$ for $\Sp$ and
$\sigma=-1$ for $\Sm$, so that the child length is $L'=2L+2\sigma$.
The child initial zero-run has length
$\lambda'=\lambda+1$ for $\Sp$ and $\lambda'=\lambda$ for $\Sm$.
The run expansions show that this is followed by one single $1$ and
at least three zeros. Indeed, the parent has an initial zero-run of
length at least four, followed by one $1$ and at least three zeros;
the next zero gap increases by one. This proves \eqref{eq:low-prefix}
for the child. Reverse-complement symmetry gives its reflected version.

It remains to exclude other low heights. Let an output prefix have
length $t$ and height $H\le\lambda'+2$.
By symmetry assume $0\le t\le L'/2$; the empty prefix is harmless.
If the zero tape has finished, the height is $L'-t$, which cannot be
this small in the chosen half. Otherwise the raw reading-head identity
in the proof of Lemma~\ref{lem:averaging} supplies parent positions $r,u$
with
\begin{equation}\label{eq:low-heads}
 H=h_{A_k}(r)+h_{A_k}(u)+\sigma,\quad
 r+u=t-\sigma,\quad r-u=H+2c-1,
\end{equation}
where $c\in\{0,1\}$ is the last output bit.
For $\Sm$ both positions are strictly internal, so both parent heights
are at least one. For $\Sp$ the ordinary parent heights are nonnegative.
Consequently, if either parent height were at least $\lambda+3$, then
$H>\lambda'+2$ in either case. Both positions must therefore lie in
the parent's two end intervals from the induction hypothesis.

They cannot be in opposite end intervals: their distance would be at
least $L-2(\lambda+4)>\lambda+4$, whereas
$|r-u|\le H+1\le\lambda+4$.
Here $L>3\lambda+12$ holds at $k=8$ and persists under the length
recursion. Nor can both positions be in the right end interval, since
$r+u\le L+1$ but then $r+u\ge2L-2\lambda-8>L+1$.
Thus both lie in $[0,\lambda+4]$.

Put $b=\mathbf1_{\{r>\lambda\}}+\mathbf1_{\{u>\lambda\}}$.
Equations~\eqref{eq:low-prefix} and \eqref{eq:low-heads} imply
\[
 H=t-2b,\qquad u=b-c+(1-\sigma)/2\le3.
\]
Since $\lambda\ge4$, the position $u$ contributes zero to $b$;
therefore $b\le1$. The assumed height bound now forces
$t\le\lambda'+4$. This excludes all low heights outside the claimed
end intervals and completes the induction.
\end{proof}

\subsection{Uniform extra-depth control and the crossing formulas}

A fixed parameter changes the height thresholds by a constant.
The low-height lemma must therefore be strengthened before the finite
crossing formulas can be transferred.

\begin{lemma}[Low heights to any fixed extra depth]\label{lem:family-low}
Let $c\ge2$ be an integer and $k\ge c+6$. Set
$\lambda=\lfloor k/2\rfloor$, $L=\ell_k$. At integer positions
$0\le t\le\lambda+c+2$,
\begin{equation}\label{eq:family-low}
 h_{A_k}(t)=
 \begin{cases}t,&t\le\lambda,\\t-2,&\lambda+1\le t\le\lambda+c+2.
 \end{cases}
\end{equation}
The reflected formula holds at the other end, and every other height
is at least $\lambda+c+1$. These end intervals therefore contain
all heights at most $\lambda+c$.
\end{lemma}
\begin{proof}
The case $c=2$ is Lemma~\ref{lem:low-heights}. For every $k\ge4$ the
first three constant runs of $A_k$ are
\[
 0^{\lfloor k/2\rfloor}\,1\,0^{k-3}.
\]
This holds for $A_4=001011$ and follows by induction from the two run
expansions: the first run changes as stated in Lemma~\ref{lem:runs},
the intervening one stays single, and the next zero-run gains one.
Thus \eqref{eq:family-low} is already exact in its stated end intervals.
It remains to exclude other low positions.

Induct on $c\ge3$, using the assertion at depth $c-1$ for the parent
word $A_{k-1}$. Write $\mu=\lfloor(k-1)/2\rfloor$ and let
$\sigma=1$ for the plus map, $-1$ for the minus map. By symmetry take
an output prefix in the left half. If the zero tape has finished,
its height is its remaining length and exceeds the proposed threshold.
Otherwise the ordinary head positions satisfy
\[
 H=h(r)+h(u)+\sigma,\qquad r+u=t-\sigma,\qquad
 r-u=H+2\chi-1,
\]
where $\chi\in\{0,1\}$ is the last bit. In the minus case both parent
positions are strictly internal. Since $H>0$ at an internal child
prefix, these identities also give $r\ge u$.

Suppose $H\le\lambda+c$. If both parent heights are at most
$\mu+c-1$, the induction hypothesis puts their positions in the
parent's two end intervals. They cannot be in opposite ends, because
their distance would exceed $H+1$, nor can both be in the right end,
because $r+u\le\ell_{k-1}+1$. The needed inequalities follow from
$\ell_{k-1}>3(\mu+c+2)$ for $k\ge c+6$, $c\ge3$.
Hence both positions are in the left end interval.

If one parent height is at least $\mu+c$, the other is zero in the
plus case or at most one in the minus case. In the plus case this
forces the smaller position $u=0$; the alternatives at the right end
contradict the left-half sum bound. In the minus case it forces $u=1$:
by the induction hypothesis the only parent positions of height one
are $1,\ell_{k-1}-1$, and again the right-end alternative is excluded.
The head-distance bound now puts both positions in
$[0,\mu+c+2]$. Their heights are given by the leading-run formula,
including its last possible extra position.

In either case let $b=\mathbf1_{\{r>\mu\}}+
\mathbf1_{\{u>\mu\}}$. The leading-run formula implies
$H=t-2b$, and the head identities give
$u=b-\chi+(1-\sigma)/2\le3$. Since $\mu\ge4$, the smaller position
cannot contribute to $b$; thus $b\le1$. It follows that
$t\le\lambda+c+2$. Reflection completes the exclusion and the induction.
\end{proof}

For an elementary factor of $U_{p,k}$, denote its beginning in the
\emph{core index} by $B$, its elementary level by $l$, and its sign
by $\varepsilon$. The complete list of factors is
\begin{equation}\label{eq:family-factor-positions}
\begin{array}{c|c|c}
 B&l&\varepsilon\\\hline
 s_{p,k}+i2^k&k&+1\quad(0\le i\le p-2)\\
 s_{p,k}+i2^k+\ell_k&k+1&-1\quad(0\le i\le p-2)\\
 s_{p,k}+(p-1)2^k&k+2&+1.
\end{array}
\end{equation}

\begin{theorem}[All finite zeros and crossings for the higher members]
\label{thm:family-zeros}
Fix $p\ge3$, put $d=\lfloor\log_2p\rfloor$, and take
$j\ge J_p:=p+d+10$. All zeros $Q_p(n)-n/2=0$ on
$2^j\le n<2^{j+1}$ are obtained from the following table, retaining
only the indices in that block. Write $D=d_{p,k}$.
\begin{center}
\begin{tabular}{@{}lll@{}}\toprule
Macro level and factor & left even zero & right even zero\\\midrule
$k$ odd, $\varepsilon=-1$
 &$2(B+D+3)$&$2(B+\ell_l-D-3)$\\
$k$ even, $\varepsilon=+1$
 &$2(B+D+3)$&$2(B+\ell_l-D-5)$\\\bottomrule
\end{tabular}
\end{center}
There are no zeros on the other factors. For odd-parity sign changes,
use positive factors at odd $k$ and negative factors at even $k$.
The left indices of the two crossing pairs $(n,n+2)$ on such a factor
are exactly
\begin{equation}\label{eq:family-odd-crossings}
 n=2(B+D+3)-1,\qquad n=2(B+\ell_l-D-4)-1.
\end{equation}
There are exactly $2p-1$ even zeros and exactly $2p-1$ sign changes
on each parity subsequence in every such standard block. Every even
zero is transverse; there are no tangencies. Thus the four/six-zero
alternation of Theorem~\ref{thm:finite-zeros} is special to $p=2$.
All these indices are within distance $2j+4p+16$ of one of the centers
$2^jz$, $z\in\mathcal Z_p$. The normalized transition-window order
$O_p(j2^{-j})$ is sharp.
\end{theorem}
\begin{proof}
The exact error identities on a signed elementary factor are
\begin{align*}
 Q_p(2m)-m&=
 \begin{cases}
 D+1+\varepsilon h_{A_l}(m-B),&k\text{ odd},\\
 D+2-\varepsilon h_{A_l}(m-B+1),&k\text{ even},
 \end{cases}\\
 Q_p(2m-1)-(2m-1)/2&=
 \begin{cases}
 D+3/2-\varepsilon h_{A_l}(m-B),&k\text{ odd},\\
 D+3/2+\varepsilon h_{A_l}(m-B),&k\text{ even}.
 \end{cases}
\end{align*}
At factor joins, heights zero or one give the consistent adjacent
core formulas; all roots below occur strictly inside a factor.
Let $\lambda=\lfloor l/2\rfloor$. On a factor selected for an even
zero, its required height exceeds $\lambda$ by $p-2$ or $p-1$.
For $p\ge3$ this is positive. Lemma~\ref{lem:family-low}, with extra
depth $c=p$, therefore gives exactly the two height positions
$T+2,\ell_l-T-2$ for the required integer height $T$.
The one-position shift at even macro levels gives the two rows of
the table. The required half-integer height for an odd crossing is
$D+3/2$; its increasing and decreasing crossing positions are
$D+3,\ell_l-D-4$, giving \eqref{eq:family-odd-crossings}.

The macro levels intersecting the standard block are
$j-d-2$ and $j-d-1$, up to their exact displaced boundary.
In particular all relevant elementary levels exceed $p+6$ under
$j\ge J_p$. The low-height lemma applies everywhere needed.
At each join exactly one adjacent factor has negative error in its
interior, for each parity. All even hits are above the small local
minimum at height $\lambda-1$, so they are transverse. Consequently
each join supplies one zero and one crossing to the respective parity
subsequence, with no additional ones in the interiors.

For an explicit separation argument, the numbers in $\mathcal Z_p$
have a common denominator $3\cdot2^{d+1}$ and none is one or two.
Thus their mutual distances and their distances to the standard-block
ends are at least $1/(3\cdot2^{d+1})$.
From \eqref{eq:family-skeleton} and the bounded corrections in $\ell_l$,
every tabulated index is within $2j+4p+16$ of its ideal join.
The inequality
\[
 2^{j-d-1}/3>2(2j+4p+16)
\]
holds at $j=J_p$ for $p\ge3$ and persists with $j$.
Thus the neighborhoods are disjoint and inside the standard block;
there are exactly $2p-1$ of them by Theorem~\ref{thm:family-profile}.
Finally a right crossing at the end of a negative elementary factor
has displacement $-2k+O_p(1)$ from its limiting endpoint.
Since $j=k+O_p(1)$, a uniform $o_p(j2^{-j})$ window is impossible.
\end{proof}

The threshold $J_p$ is deliberately sufficient, not asserted to be
minimal. The displayed coordinate tables are exact formulas in powers,
floors, and $p$; they require neither numerical root searches nor a
quantitative Gaussian-convergence rate.

\section{Ordered frequencies and parameter-dependent dyadic laws}\label{sec:family-frequencies}
Let $\nu_p(v)=\#\{n\ge1:Q_p(n)=v\}$. These are multiplicities
in the complete infinite sequence, finite by
Theorem~\ref{thm:family-global}. A natural value block is
$\{2u+1:a_{p,k}\le u<a_{p,k+1}\}$.
We first retain the order of its frequencies; a histogram alone would
not determine a cut through a dyadic value interval.

\subsection{Plateaus and elementary frequency words}
Write $L_{F_p}(u)=\#\{m:F_p(m)=u\}$ and define $L_{G_p}$ likewise.
Then $\nu_p(2u+1)=L_{F_p}(u)+L_{G_p}(u)$.
For $U_{p,k}$ let $\alpha_i,\beta_i$ be the gap lists defined in
Section~\ref{sec:words}, and put $\beta_0=0$.
The outer runs of $U_{p,k}$, for $k\ge2$, have lengths
$\lfloor k/2\rfloor$ and $\lfloor k/2\rfloor+1$, independently of $p$.
This follows from $U_{p,2}=(0110)^{p-1}001011$ and
Lemma~\ref{lem:runs}.

\begin{lemma}[Interior and boundary local times]\label{lem:plateau}
For every $p\ge2$, $k\ge2$, and $2\le i\le L_{p,k}/2$,
\begin{equation}\label{eq:interior-frequency}
 \nu_p(2(a_{p,k}+i-1)+1)=2+\alpha_i+\beta_{i-1}.
\end{equation}
At the initial level, $\nu_p(2a_{p,k}+1)=k+1$.
\end{lemma}
\begin{proof}
On an odd block $F_p$ counts ones and $G_p$ counts zeros. An interior
plateau contributes one plus the corresponding zero gap and one plus
the preceding one gap. On an even block the roles exchange; the
one-position shift in $G_p$ leaves complete interior plateau lengths
unchanged. Only the initial level can extend into the preceding block.
For $k=2r+1$ its $F_p$ plateau has $r+1$ points on each side, with the
shared endpoint counted once, hence length $k$; the other plateau has
length one. For $k=2r$ the $G_p$ plateau has $r+1$ and $r$ points,
again with one overlap, hence length $k$; its partner has length one.
The case $k=2$ uses the terminal one-run of length one in $U_{p,1}$.
This proves the boundary formula. The clock offset $p-2$ changes no
plateau length.
\end{proof}

\begin{lemma}[Disjoint gap supports]\label{lem:gap-support}
For every elementary word $A_k$, $k\ge2$, exactly one of $\alpha_i,\beta_{i-1}$ is positive
at each $1\le i\le\ell_k/2$.
\end{lemma}
\begin{proof}
We prove the property for an output of either map; the words $A_2=A_3=01$
verify the initial cases.
For a parent $P$ of length $2p$, let $\operatorname{pos}(0_j)$ and
$\operatorname{pos}(1_j)$ be its zero-based symbol positions.
In $\Sp(P)$, the positive $\alpha$ gaps occur at index $1$ and at
$2+\operatorname{pos}(0_j)$, while the positive $\beta_{i-1}$ gaps occur
at $2+\operatorname{pos}(1_j)$.
The final output one-run gives an additional index $2p+2$, outside
the range $1\le i\le2p+1$. The listed supports partition that range.
For $\Sm(P)$ the two supports are
$1+\operatorname{pos}(0_j)$ and $1+\operatorname{pos}(1_j)$, with the
terminal index $2p$ omitted from $1\le i\le2p-1$.
All indicated run lengths are positive by
\eqref{eq:plus-runs}--\eqref{eq:minus-runs}.
The exceptional parent $01$ is checked directly.
\end{proof}

For the elementary word $A_k$ write
\[
 q_k=\ell_k/2=\frac{2^{k-1}+(-1)^k}{3},\qquad
 \lambda_k=\lfloor k/2\rfloor.
\]
This elementary semilength $q_k$ is independent of the family parameter. Let $\alpha_i^{(k)},\beta_i^{(k)}$ be the two gap lists of $A_k$,
with $\beta_0^{(k)}=0$. Define an integer word $\mathbf V_k$ of length
$q_k$ by
\begin{equation}\label{eq:V-frequency-word}
 \mathbf V_k[1]=k+1,\qquad
 \mathbf V_k[i]=2+\alpha_i^{(k)}+\beta_{i-1}^{(k)}
 \quad(2\le i\le q_k).
\end{equation}
These integer words are indexed from one, unlike the binary words.
The number $q_k$ is odd. Consequently the symmetry below lets us write
\begin{equation}\label{eq:V-palindrome}
 \mathbf V_k=(k+1)\,\mathbf R_k\,\operatorname{rev}(\mathbf R_k),
 \qquad |\mathbf R_k|=(q_k-1)/2.
\end{equation}
The operation $\operatorname{rev}$ reverses an integer word; it is not
bitwise complementation. A half-word may be empty.

\begin{lemma}[Elementary frequency symmetry and polynomial]\label{lem:frequency-words}
Every $\mathbf V_k$ has the form \eqref{eq:V-palindrome}. Define
\[
 \mathcal B_k(z)=\sum_{i=1}^{q_k}z^{\mathbf V_k[i]}.
\]
Then
\begin{equation}\label{eq:elementary-frequency-rec}
 \mathcal B_2(z)=z^3,\qquad
 \mathcal B_{k+1}(z)=zB_k(z)+(q_k-(-1)^k)z^3\quad(k\ge2).
\end{equation}
Its leading coefficient is one, at degree $k+1$, and its coefficient
at degree $k$ is zero.
\end{lemma}
\begin{proof}
Reverse-complement symmetry of $A_k$ gives
$\beta_{i-1}^{(k)}=\alpha_{q_k+2-i}^{(k)}$ for $2\le i\le q_k$.
Thus the entries with indices $2,\ldots,q_k$ form a palindrome of even
length, proving \eqref{eq:V-palindrome}.

For the polynomial statement let $T_k^\circ(z)$ enumerate all constant
runs of $A_k$, counting both symbols. The run expansions and the
disjoint-support proof of Lemma~\ref{lem:gap-support}, applied now to
$A_k$, show that
\begin{equation}\label{eq:E-run-polynomial}
 \mathcal B_k(z)=z^2T_k^\circ(z)-2z^{\lambda_k+2}+z^{k+1}.
\end{equation}
The two subtracted terms are the omitted final one-run and the first
zero-run, whose frequency is replaced by $k+1$.
The initial words $A_2=A_3=01$ satisfy the same identity directly.
The number of constant runs is $q_k+1$. The run expansions give
\[
 T_{k+1}^\circ(z)=
 \begin{cases}
 zT_k^\circ(z)+(q_k+1)z,&k\text{ odd},\\
 zT_k^\circ(z)+(q_k-1)z+2(1-z)z^{\lambda_k},&k\text{ even}.
 \end{cases}
\]
For $A_2=01$ the second formula is checked separately; it remains valid.
Substitution in \eqref{eq:E-run-polynomial} proves
\eqref{eq:elementary-frequency-rec}. The leading coefficient and the
missing next coefficient follow by induction, beginning with
$\mathcal B_2=z^3$ and $\mathcal B_3=z^4$.
\end{proof}

\subsection{Ordered natural frequencies: the repeated-join correction}

Use the elementary integer
frequency words $\mathbf V_k$, their palindromic tails $\mathbf R_k$,
and polynomials $\mathcal B_k(z)$ of Lemma~\ref{lem:frequency-words}.
Let $\widetilde{\mathbf V}_k$ be obtained by increasing only the first
entry of $\mathbf V_k$ from $k+1$ to $k+2$.
This distinction is essential at the new repeated-pair joins.

\begin{theorem}[Ordered family frequencies]\label{thm:family-natural-frequency}
For every $p\ge2$ and $k\ge2$, the ordered frequencies on
$a_{p,k}\le u<a_{p,k+1}$ are exactly
\begin{equation}\label{eq:family-frequency-word}
 \mathbf V_k\mathbf V_{k+1}
 (\widetilde{\mathbf V}_k\mathbf V_{k+1})^{p-2}\mathbf V_{k+2}.
\end{equation}
Consequently their polynomial is
\begin{equation}\label{eq:family-natural-polynomial}
 P_{p,k}(z)=(p-1)(\mathcal B_k(z)+\mathcal B_{k+1}(z))+\mathcal B_{k+2}(z)
 +(p-2)(z^{k+2}-z^{k+1}).
\end{equation}
If $v_{p,k}=a_{p,k+1}-a_{p,k}=L_{p,k}/2$, then
\begin{equation}\label{eq:family-natural-rec}
 \begin{split}
 P_{p,2}(z)&=3z^3+(2p-3)z^4+z^5,\\
 P_{p,k+1}(z)&=zP_{p,k}(z)+(v_{p,k}-(-1)^k)z^3
 \qquad(k\ge2).
 \end{split}
\end{equation}
The separate initial block has $P_{p,1}(z)=z^3+(p-1)z^4$; its
step to $k=2$ must not be inferred from the displayed recursion.
For every $k\ge1$ the total multiplicity on that block is
\[
 4v_{p,k}-\frac{3+(-1)^k}{2}.
\]
For $k\ge2$ the largest multiplicity is $k+3$, attained once, and
\[
 \frac{P_{p,k}(z)}{v_{p,k}}\longrightarrow\frac{z^3}{2-z}
 \quad\text{in coefficient }\ell^1\text{ norm, with error }O_p(k2^{-k}).
\]
In particular the limiting mean and variance are four and two.
\end{theorem}
\begin{proof}
By Lemma~\ref{lem:plateau}, a
complete interior plateau has frequency $2+\alpha_i+\beta_{i-1}$.
The constant $p-2$ in the clock does not change plateau lengths.
The first frequency at a macro boundary is still $k+1$, because the
initial and final constant runs of $U_{p,k}$ have lengths
$\lfloor k/2\rfloor$ and $\lfloor k/2\rfloor+1$ independently of $p$.
Their two-block overlap is included in Lemma~\ref{lem:plateau},
also for $k=2$.

At a join from $A_k$ to $\overline{A}_{k+1}$ the frequency of the
joining level is $2+\lfloor k/2\rfloor+\lfloor(k+1)/2\rfloor=k+2$.
The reverse join to another $A_k$ gives the same $k+2$, not $k+1$;
this is the changed first entry in $\widetilde{\mathbf V}_k$.
At the final join to $A_{k+2}$ it is $k+3$.
All remaining entries are the corresponding unchanged elementary
frequencies, also on a complemented factor because its two gap
contributions exchange without reversing value order.
This proves \eqref{eq:family-frequency-word} and its polynomial.

Insert the scalar elementary recurrence \eqref{eq:elementary-frequency-rec} in
\eqref{eq:family-natural-polynomial}. The additional monomials advance
by multiplication by $z$, and the new $z^3$ coefficient is
$v_{p,k}-(-1)^k$. Directly $\mathcal B_2=z^3$, $\mathcal B_3=z^4$,
$\mathcal B_4=z^5+2z^3$, giving the stated initial polynomial.
Reading the first word separately gives $P_{p,1}$.
Differentiation of \eqref{eq:family-natural-rec}, starting at $k=2$,
proves the mass formula; $k=1$ is immediate from its explicit polynomial.
The maximal degree and its coefficient follow from
\eqref{eq:family-natural-polynomial}: only $\mathcal B_{k+2}$ has degree $k+3$.
For completeness, put $G(z)=z^3/(2-z)$ and
$E_k=P_{p,k}-v_{p,k}G$. Since $v_{p,k+1}=2v_{p,k}-(-1)^k$,
\[
 E_{k+1}=zE_k+(-1)^k(G-z^3).
\]
The coefficient norm of $G-z^3$ is one, so $\|E_k\|_1=O_p(k)$.
Division by $v_{p,k}\asymp_p2^k$ proves the stated rate. To pass to
the first two moments, on degrees at most $k+3$ multiply the
coefficient error by at most $(k+3)^2$; the remaining geometric tail
has second moment $O(k^2 2^{-k})$. Both moments therefore converge to
$G'(1)=4$ and $G''(1)+G'(1)-G'(1)^2=2$.
\end{proof}

\subsection{Worked example: \texorpdfstring{$p=3$}{p=3}, a repeated join, and a dyadic cut}
This example separates three operations: forming a signed binary word,
recording its ordered multiplicities, and cutting at a power of two.
For $p=3$ and elementary level $k=4$,
\[
 U_{3,4}=A_4\overline{A_5}A_4\overline{A_5}A_6,
 \qquad (\ell_4,\ell_5,\ell_6)=(6,10,22).
\]
Its natural value-index block is $26\le u<53$, since
$a_{3,4}=26$ and $a_{3,5}=53$. The elementary frequency words are
\[
\begin{split}
 \mathbf V_4&=(5,3,3),\qquad
 \mathbf V_5=(6,4,3,3,4),\\
 \mathbf V_6&=(7,5,3,3,4,3,3,4,3,3,5).
\end{split}
\]
The first $A_4$ contributes $\mathbf V_4$, but the second contributes
$\widetilde{\mathbf V}_4=(6,3,3)$, not $(5,3,3)$.
Indeed its preceding factor is $\overline{A_5}$, and the two runs at
this reverse join have lengths $\lambda_5=2$ and $\lambda_4=2$.
The joining multiplicity is $2+2+2=6$ rather than the distinguished
first entry $5$ of $\mathbf V_4$. The complete natural frequency word
is therefore
\[
 \mathbf V_4\,\mathbf V_5\,
 \underbrace{(6,3,3)}_{\widetilde{\mathbf V}_4}\,
 \mathbf V_5\,\mathbf V_6.
\]
The five factor beginnings are $u=26,29,34,37,42$.
Their first frequencies are $5,6,6,6,7$, respectively. In particular
the extra join changes the multiplicity of the actual value
$2\cdot34+1=69$ from the uncorrected prediction $5$ to $6$.
Taking a histogram gives
\[
 P_{3,4}(z)=14z^3+6z^4+3z^5+3z^6+z^7.
\]

Now cut at the pure dyadic endpoints $32$ and $64$, so $j=5$.
Here $d=1$, $b=2$, $k=4$, and the two integers in the nonpower case
are $t=2b-p=1$ and $u=2p-2b-1=1$.
The lower cut lies in the first $\mathbf V_5$, starting at $29$:
its distance is $32-29=3=(q_5+1)/2$.
As $\mathbf V_5=6\,(4,3)\,(3,4)$, the retained suffix is
$\operatorname{rev}(\mathbf R_5)=(3,4)$.
In the next natural block, $\mathbf V_6$ starts at $58$, and
$64-58=6=(q_6+1)/2$; thus the retained prefix is
$7\mathbf R_6=(7,5,3,3,4,3)$.
Equation~\eqref{eq:family-pure-ordered-nonpower} becomes the explicit
concatenation
\[
 \underbrace{(3,4)}_{\operatorname{rev}(\mathbf R_5)}\,
 \underbrace{\widetilde{\mathbf V}_4\mathbf V_5}_{8\text{ entries}}\,
 \mathbf V_6\,\mathbf V_5\,
 \underbrace{(7,5,3,3,4,3)}_{7\mathbf R_6}.
\]
Its $32$ entries have histogram
\[
 D_{3,5}(z)=16z^3+8z^4+3z^5+3z^6+2z^7.
\]
The two maximal multiplicities, both seven, occur at $u=42,58$;
the three second-largest multiplicities, all six, occur at
$u=34,37,53$. These positions distinguish the ordered result from
its histogram. The mean is $4-2^{-5}$, the variance is $1567/1024$,
and the total-variation distance from $z^3/(2-z)$ is $1/16$,
as given by Corollary~\ref{cor:family-frequency-statistics}.

\subsection{Closed formulas on pure dyadic value blocks}

Define
\[
 D_{p,j}(z)=\sum_{u=2^j}^{2^{j+1}-1}z^{\nu_p(2u+1)},\qquad
 d=\lfloor\log_2p\rfloor,\quad b=2^d,\quad k=j-d.
\]
The following formulas cover every sufficiently large dyadic value block
for each $p$, with an explicit threshold $j\ge d+2$.
The smaller, finitely many blocks are fixed by the initial construction.

\begin{theorem}[Parameter-dependent pure dyadic histograms]
\label{thm:family-dyadic}
For $j\ge d+2$, if $p=b$ is a power of two, then
\begin{equation}\label{eq:family-dyadic-power}
 D_{p,j}(z)=\sum_{r=3}^{k+1}2^{j+2-r}z^r
       +2z^{k+2}+(2p-3)z^{k+3}+z^{k+4}.
\end{equation}
If $b<p<2b$, then
\begin{equation}\label{eq:family-dyadic-nonpower}
\begin{split}
 D_{p,j}(z)={}&\sum_{r=3}^{k}2^{j+2-r}z^r
  +(4b-2p+1)z^{k+1}\\
 &+(4p-4b-1)z^{k+2}+(4b-2p)z^{k+3}.
\end{split}
\end{equation}
An empty sum is zero. All displayed coefficients are positive integers.
For each fixed $p$, the normalized law tends to $G(z)=z^3/(2-z)$.
For every $p>2$, however, its finite histogram differs from
$\{3+\nu_2(a):1\le a\le2^j\}$ on all blocks in the stated range.
\end{theorem}
\begin{proof}
The lower cut $2^j$ lies in macro block $k=j-d$, since
$c_p<6b<2c_p$. There are two possible geometries.

If $p=b$, the cut bisects the final word $\mathbf V_{k+2}$, in the
sense that it follows its distinguished first entry and its first
palindromic half. The retained suffix is $\operatorname{rev}(\mathbf R_{k+2})$.
The upper cut is in the corresponding factor of the next macro block,
retaining its first entry and first half. Thus the full ordered word is
\begin{equation}\label{eq:family-pure-ordered-power}
\begin{split}
 \operatorname{rev}(\mathbf R_{k+2})\,
 \mathbf V_{k+1}\mathbf V_{k+2}
 (\widetilde{\mathbf V}_{k+1}\mathbf V_{k+2})^{p-2}
 (k+4)\mathbf R_{k+3}.
\end{split}
\end{equation}
This includes exactly the three-word case at $p=2$.

If $p>b$, put $t=2b-p\ge1$, $u=2p-2b-1\ge1$.
The lower cut is in the $t$-th complemented factor, of level $k+1$.
Its exact starting value index is
\[
 a_{p,k}+(t-1)2^{k-1}+\ell_k/2
 =\frac{(6b-1)2^k-3+(-1)^k}{6}.
\]
Its distance to $2^j$ is $(q_{k+1}+1)/2$, where $q_l=\ell_l/2$.
This is precisely the same palindromic half-cut, not an asymptotic cut.
The whole ordered word is
\begin{equation}\label{eq:family-pure-ordered-nonpower}
\begin{split}
 &\operatorname{rev}(\mathbf R_{k+1})
 (\widetilde{\mathbf V}_k\mathbf V_{k+1})^u\mathbf V_{k+2}\\
 &\qquad\mathbf V_{k+1}
 (\mathbf V_{k+2}\widetilde{\mathbf V}_{k+1})^{t-1}
 (k+3)\mathbf R_{k+2}.
\end{split}
\end{equation}
For the power case the analogous distance is $(q_{k+2}+1)/2$.
These exact distances prove both ordered statements, including endpoint
placement and every distinguished-entry correction.

The polynomial of $\mathbf R_l$ is $(\mathcal B_l-z^{l+1})/2$. Hence
\eqref{eq:family-pure-ordered-power} gives
\begin{equation}\label{eq:family-pure-B-power}
\begin{split}
 D_{p,j}={}&(p-1)\mathcal B_{k+1}+(p-\tfrac12)\mathcal B_{k+2}+\tfrac12 \mathcal B_{k+3}\\
 &-(p-2)z^{k+2}+(p-\tfrac52)z^{k+3}+\tfrac12z^{k+4},
\end{split}
\end{equation}
whereas \eqref{eq:family-pure-ordered-nonpower} gives
\begin{equation}\label{eq:family-pure-B-nonpower}
\begin{split}
 D_{p,j}={}&uB_k+(p-\tfrac12)\mathcal B_{k+1}+(t+\tfrac12)\mathcal B_{k+2}\\
 &-uz^{k+1}+(u-t+\tfrac12)z^{k+2}+(t-\tfrac12)z^{k+3}.
\end{split}
\end{equation}
One may simplify them using the exact rational expression
\begin{equation}\label{eq:elementary-B-closed}
 \mathcal B_l(z)=\frac{2^l z^3}{6(2-z)}
       +\frac{2(-1)^l z^3}{3(z+1)}
       +\frac{z^{l+2}(z-1)}{(z-2)(z+1)}\qquad(l\ge2).
\end{equation}
This identity follows from $\mathcal B_2=z^3$ and \eqref{eq:elementary-frequency-rec}; its apparent poles
cancel because the left side is a polynomial.
In both combinations the alternating term cancels. The coefficients
remaining at the last three degrees are precisely those in
\eqref{eq:family-dyadic-power} or
\eqref{eq:family-dyadic-nonpower}; all lower coefficients are
$2^{j+2-r}$. This proves the formulas.
Their coefficients sum to $2^j$, and their tail discrepancy from
$2^jG$ has bounded coefficient norm for fixed $p$, proving convergence.
For $p>2$ their maximal degree is strictly below $j+3$, so the original
shifted-valuation histogram is not equal to them.
\end{proof}

For example, in their respective stated ranges,
\begin{align}
 D_{3,j}&=D_{2,j}-z^j+z^{j+1}+z^{j+2}-z^{j+3},\\
 D_{4,j}&=D_{2,j}-2z^j+3z^{j+1}-z^{j+3}.
\end{align}
Thus $Q_3$ and $Q_4$ preserve every sufficiently low multiplicity count
exactly, but redistribute the top of the histogram. Treating them as
unchanged copies of Theorem~\ref{thm:power-blocks} would be incorrect.

\begin{corollary}[Exact moments, distance, and cumulative count]
\label{cor:family-frequency-statistics}
Choose $u$ uniformly from $[2^j,2^{j+1})$, with $j\ge d+2$, and set
$X_{p,j}=\nu_p(2u+1)$. Then
\begin{align}
 \mathbb E X_{p,j}&=4-2^{-j},\\
 \operatorname{Var}(X_{p,j})&=
 2-\bigl(2(j-d)+4p-5\bigr)2^{-j}-4^{-j}.
 \label{eq:family-variance}
\end{align}
The total-variation distance to the geometric limit is exactly
\begin{equation}\label{eq:family-TV}
 d_{\rm TV}(X_{p,j},G)=T_p2^{-j},\quad
 T_p=
 \begin{cases}
 \frac{3p}{2}-3+(1-p/4)_+,&p=b,\\
 (4p-5b-1)_+ +(7b/2-2p)_+,&b<p<2b,
 \end{cases}
\end{equation}
where $x_+=\max(x,0)$. The coefficient $\ell^1$ distance is twice
this number. Also the count over the entire infinite sequence is
\begin{equation}\label{eq:family-cumulative-frequency}
 \#\{n\ge1:Q_p(n)<2^{j+1}\}
 =2^{j+2}-j+d+2-2p.
\end{equation}
This does not assert that the counted indices form an initial interval.
\end{corollary}
\begin{proof}
Differentiation of the finite polynomials yields
$\mathbb EX=4-2^{-j}$ and
$\mathbb EX^2=18-(2k+4p+3)2^{-j}$, proving the variance.
For completeness, in the nonpower case the nonzero leading
coefficients of $D_{p,j}-2^jG$, at degrees $k+1,k+2,k+3$, are
\[
 -(2p-2b-1),\quad 4p-5b-1,\quad 7b/2-2p.
\]
All remaining coefficients are a negative geometric tail. In the
power case the corresponding coefficients at $k+2,k+3,k+4$ are
\[
 2-p,\quad3p/2-3,\quad1-p/4,
\]
again followed by a negative tail. Their positive mass divided by
$2^j$ gives \eqref{eq:family-TV} exactly.

For the cumulative formula, the initial values $u=0,\ldots,p-2$
have multiplicity two and $u=p-1$ has multiplicity four. Up to the
beginning $a_{p,k}$, sum these $2p+2$ occurrences and the complete
natural-block masses in Theorem~\ref{thm:family-natural-frequency}.
At the lower pure-dyadic cut, add the complete elementary words and
the first palindromic half identified in the preceding proof.
The required elementary masses are
\[
 \mathcal B_l'(1)=4q_l-\frac{1+(-1)^l}{2},\qquad
 \text{first-half mass}=\frac{\mathcal B_l'(1)+l+1}{2};
\]
each tilde word adds one. In the power case the partial contribution
is
\[
 (p-1)(\mathcal B_k'(1)+\mathcal B_{k+1}'(1))+(p-2)
   +\frac{\mathcal B_{k+2}'(1)+k+3}{2},
\]
and in the nonpower case it is
\[
 tB_k'(1)+(t-1)\mathcal B_{k+1}'(1)+(t-1)
   +\frac{\mathcal B_{k+1}'(1)+k+2}{2}.
\]
Substitution of $q_l$, $a_{p,k}$ and $j=k+d$ in either expression,
including the complete preceding masses, gives
$4\cdot2^j-k+2-2p$, which is \eqref{eq:family-cumulative-frequency}.
\end{proof}

\begin{proposition}[All locations of the two largest multiplicities]
\label{prop:family-extrema}
Continue to assume $j\ge d+2$, $k=j-d$. If $p=b$, the maximum is
$k+4$, attained only at
\[
 u=a_{p,k+1}+(p-1)2^k.
\]
The next largest multiplicity $k+3$ occurs at exactly
\[
\begin{split}
 &a_{p,k+1}+i2^k\quad(1\le i\le p-2),\\
 &a_{p,k+1}+i2^k+q_{k+1}\quad(0\le i\le p-2).
\end{split}
\]
Its count is $2p-3$.
If $p>b$, put $t=2b-p$. The maximum $k+3$ occurs exactly at
\[
\begin{split}
 &a_{p,k}+(p-1)2^{k-1},\\
 &a_{p,k+1}+i2^k\quad(1\le i\le t-1),\\
 &a_{p,k+1}+i2^k+q_{k+1}\quad(0\le i\le t-1).
\end{split}
\]
There are $2t=4b-2p$ maximizers. The next largest multiplicity $k+2$
occurs at $a_{p,k+1}$ and at
\[
 a_{p,k}+i2^{k-1},\quad
 a_{p,k}+i2^{k-1}+q_k\qquad(t\le i\le p-2).
\]
Its count is $4p-4b-1$. Every position listed is a value index $u$;
the corresponding actual value is $2u+1$.
\end{proposition}
\begin{proof}
The only entry $l+1$ of $\mathbf V_l$ is its first; none of its
remaining entries equals $l$. In a tilde word the first entry is
$l+2$, while the remaining entries stay unchanged. Read these
positions in the exact ordered words
\eqref{eq:family-pure-ordered-power} and
\eqref{eq:family-pure-ordered-nonpower} and translate each factor
beginning to its value coordinate. The displayed lists exhaust
all occurrences; their cardinalities agree with the top two
coefficients of the finite histogram formulas.
\end{proof}

\begin{figure}[htbp]
\centering
\includegraphics[width=.99\textwidth]{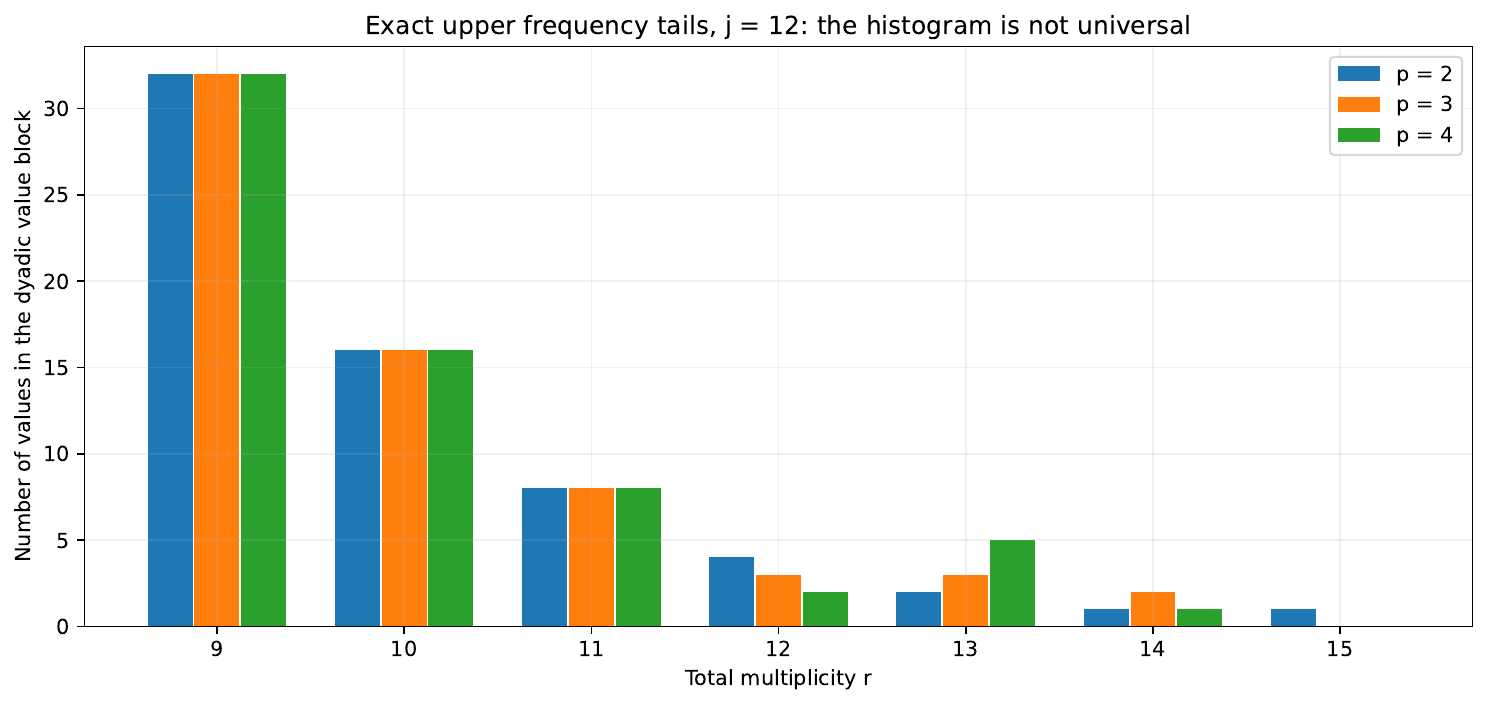}
\caption{Exact top multiplicities on the same pure dyadic value block
$j=12$, for $p=2,3,4$. The low part of the histogram agrees; the
largest multiplicities do not. In particular the $Q_3$ maximum is
attained twice, whereas the $Q_4$ maximum is attained once.}
\label{fig:family-frequency}
\end{figure}

\section{Exact summation and the periodic linear term}\label{sec:family-summation}

Put $S_p(N)=\sum_{n=1}^NQ_p(n)$ and
$\mathcal E_p(N)=S_p(N)-N(N+1)/4$.

\begin{theorem}[Parameter-dependent exact summation]\label{thm:family-sum}
Let $M\ge p+2$ and choose $K$ by $s_{p,K}\le M<s_{p,K+1}$.
With $T=M-s_{p,K}+1$,
\begin{equation}\label{eq:family-exact-sum}
\begin{split}
\mathcal E_p(2M)={}&
 (K+2p-\tfrac52)(M+1)-\sum_{k=1}^K s_{p,k}
 +\tfrac32-p^2\\
 &-\left\lfloor\frac{K-1}{2}\right\rfloor
 -\mathbf1_{\{K\text{ even}\}}h_{U_{p,K}}(T),
\end{split}
\end{equation}
where
\begin{equation}\label{eq:family-boundary-sum}
 \sum_{k=1}^K s_{p,k}
 =\frac{c_p}{3}(2^{K+1}-2)-(p-1)K
  -\frac{(-1)^K-1}{6}-\left\lfloor\frac{(K-1)^2}{4}\right\rfloor.
\end{equation}
Define $\Psi_p$ by continuous one-periodic extension of
\begin{equation}\label{eq:family-Psi}
 \Psi_p(\theta)=p-1-\frac12\log_2 C_p-\frac\theta2-2^{-\theta}
 \qquad(0\le\theta\le1).
\end{equation}
Then for all integer $N\to\infty$,
\begin{equation}\label{eq:family-sum-asym}
 S_p(N)=\frac{N^2}{4}+\frac N2\log_2N
 +N\Psi_p\!\left(\{\log_2(N/C_p)\}\right)
 +O_p\!\left(\frac N{\sqrt{\log N}}\right).
\end{equation}
At even cutoffs, retaining the exact last height refines this to
\begin{equation}\label{eq:family-sum-refined}
\begin{split}
 S_p(N)={}&\frac{N^2}{4}+\frac N2\log_2N
 +N\Psi_p\!\left(\{\log_2(N/C_p)\}\right)\\
 &-\mathbf1_{\{K\text{ even}\}}h_{U_{p,K}}(T)+O_p((\log N)^2).
\end{split}
\end{equation}
Thus the coefficient of $N\log_2N$ remains $1/2$, independently of
$p$, and the periodic shape is universal up to its explicit vertical
constant and phase scale. The improved $O_p((\log N)^2)$ remainder
on odd macro levels survives for every parameter.
\end{theorem}
\begin{proof}
The symmetric mode is
\[
 F_p(m)+G_p(m)-m=d_{p,k}+\mathbf1_{\{k\text{ even}\}}U_{p,k}[t]
 \quad(m=s_{p,k}+t).
\]
Consequently the contribution of a parity pair to $\mathcal E_p$ is
$k+2p-5/2$ on odd macro levels and
$k+2p-5/2-(h(t+1)-h(t))$ on even ones.
The first $p+1$ pairs have total error
\[
 \sum_{m=1}^p(2m-3/2)+(2p-7/2)
 =p^2+\frac{3p}{2}-\frac72.
\]
An entire even macro block subtracts one height unit, since its
last retained prefix has height one. Sum the half-open macro
blocks and telescope their clock differences. The resulting
constant is
\[
 p^2+3p/2-7/2-(2p-5/2)(p+2)=3/2-p^2,
\]
which proves \eqref{eq:family-exact-sum}. Summing
\eqref{eq:family-skeleton} proves \eqref{eq:family-boundary-sum}.

Uniformly on the last macro block,
$|h_{U_{p,K}}|=O_p(2^K/\sqrt K)$. Thus for $N=2M$,
\[
 \mathcal E_p(2M)=M(K+2p-5/2)-\frac{c_p}{3}2^{K+1}
 +O_p(K^2+2^K/\sqrt K).
\]
The actual clock differs from $(c_p/3)2^K$ by $O_p(K)$.
Replacing its level by $\lfloor\log_2(N/C_p)\rfloor$ changes this
main expression by only $O_p(K)$ in the possible short displaced
boundary intervals. Set $\theta=\{\log_2(N/C_p)\}$ and use
$(c_p/3)2^{K_0+1}=N2^{-\theta}$.
Adding $N(N+1)/4$ gives \eqref{eq:family-sum-asym}; keeping the
height gives \eqref{eq:family-sum-refined}.
The periodic function has equal endpoint values and is Lipschitz.
For odd $N$, add its last term to the formula at $N-1$ and use the
pointwise growth bound. The leading increment is $N/2+O_p(\log N)$,
and the periodic term changes by $O_p(1)$; the stated error follows.
\end{proof}

\begin{figure}[htbp]
\centering
\includegraphics[width=.99\textwidth]{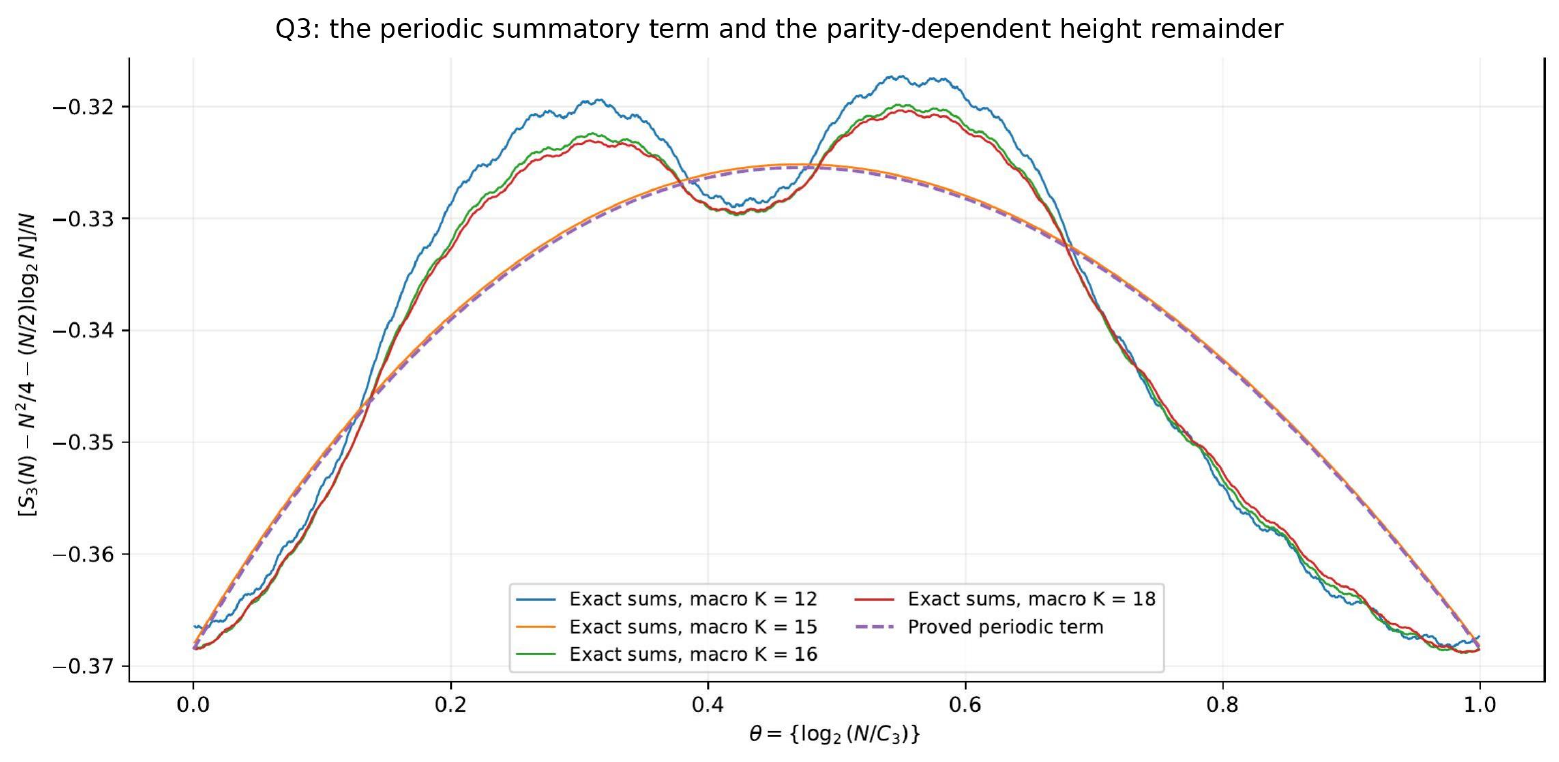}
\caption{The $Q_3$ summatory residual compared with $\Psi_3$. Odd macro
levels lose the last-height term entirely; even levels retain it.
These data use exact cumulative sums of the directly computed sequence.}
\label{fig:family-sum}
\end{figure}

\section{Comparison with the original perturbation}\label{sec:comparison}

The elementary words are related to those in Cloitre's analysis of
A394051 \cite[Section 2.3 and Lemma 4.3]{cloitre}. In that notation the positive step
word starts with $P_0=001011$, and the two word laws are
\[
N_r=\I_1(P_r[2:],P_r[:-1]),\qquad
P_{r+1}=\I_0(00N_r1,0N_r).
\]
For words in this family the designated tapes do not require a fallback.
The family construction and proofs do not assume the existence of
the original sequence or an asymptotic theorem about it.

The identification is
\begin{equation}\label{eq:original-words}
A_{2r+4}=P_r,\qquad A_{2r+5}=0N_r\qquad(r\ge0).
\end{equation}
Indeed $A_4=P_0$. When $P$ begins with two zeros, deleting the first
output zero from $\Sm(P)$ gives
$\I_1(P[2:],P[:-1])$: after its first zero the latter machine has the
same state and remaining tapes as the former after its first two zeros.
All even $A_k$, $k\ge4$, begin with at least two zeros.
Also $\Sp(0N)=\I_0(00N1,0N)$.
Induction proves \eqref{eq:original-words} from the cited word laws.

Equation~\eqref{eq:original-words} identifies the inherited elementary
excursions exactly. The exact-fit interleaving, tape counters,
height additivity, and parked-head mechanism are likewise tools from
Cloitre's work \cite{cloitre}. The structural result here is the full orbit
$U_{p,k}=(A_k\overline{A}_{k+1})^{p-1}A_{k+2}$ together with the
parameter-dependent clock $s_{p,k}$ and core reconstruction in
Section~\ref{sec:family}. The proof checks the prescribed nested recurrence
at every index. This does not follow merely from the presence of the same
elementary excursions in the original model.

The analytic method also has a specific antecedent. Uniform profiles
for the delayed perturbation were proved using nonlinear averaging
and smooth comparison functions \cite{mantovanellidelayed}.
Here Lemma~\ref{lem:averaging} supplies an exact shifted identity for
the alternating maps of lengths $2L+2$ and $2L-2$, and
Lemma~\ref{lem:comparison} proves a uniform comparison through their
boundary cases. The resulting Gaussian limit, combined with the new
block placement, gives every profile \eqref{eq:family-standard-profile} and its sharp
constant $\Lambda_p$. The four-piece profile \eqref{eq:profile}
is the specialization at $p=2$.
The exact finite zeros and crossings in Theorem~\ref{thm:family-zeros}
and the exceptional case in Theorem~\ref{thm:finite-zeros}
require the separate low-height analysis and are not consequences of
uniform convergence alone.

For $p=2$, the pure dyadic frequency histogram in
Theorem~\ref{thm:power-blocks}
is exactly the same shifted-valuation multiset as in the original
perturbation's frequency law \cite{mantovanelli}, using each model's
own value-index convention. This does not identify their order:
the $Q_2$ block $u=4,5,6,7$ has frequencies $(3,4,5,3)$, whereas
the original $s=4,5,6,7$ block, indexed by values $2s-1$, has
frequencies $(3,3,5,4)$. The exact equality of histograms is therefore
compatible with different ordered frequency words and different
arch arrangements. The ordered factorization \eqref{eq:family-frequency-word}, specialized
in \eqref{eq:ordered-natural-frequencies}, retains the boundary corrections
and the exact cuts through the palindromic tails.
It proves the valuation histogram for $Q_2$ without transferring it
from the original model. For every $p>2$, Theorem~\ref{thm:family-dyadic}
shows explicitly how the upper finite histogram differs, despite the
common limiting law.

Finally, logarithmic summatory bias and periodic linear terms already
occur for the original perturbation \cite{mantovanellizeta}.
Section~\ref{sec:family-summation} supplies the exact telescoping identity
for every clock $s_{p,k}$, the explicit function $\Psi_p$, and the
parity-dependent height remainder \eqref{eq:family-sum-refined}. These formulas distinguish the
two models even where their frequency histograms agree.

\section{Reproducibility and further questions}\label{sec:audit}
The reproducibility archive \cite{zenodoarchive} accompanies the full
family, including the three-arch specialization. It contains independent
direct and word-based generators, the exact formulas, machine-readable
audit reports, and the data and scripts for all figures. The principal
family implementation is \texttt{family.py}; \texttt{audit\_family.py}
regenerates \texttt{report\_family\_verification.json}.
The earlier three-arch verification programs and reports are retained
with their own documentation. The archive's audit guide records each
program's scope, parameter ranges, tested endpoints, and report provenance,
rather than repeating the full program catalogue here.

The family report compares every term of the two constructions for
$23$ parameter choices: $p=2,\ldots,9$ through $2^{23}$ terms each, and
$15$ additional parameters through $2^{19}$ terms each. It checks the
word factorization, range and two-step bounds, complete natural and
pure dyadic frequency blocks, extreme locations, and every even summatory
cutoff in those runs. Transition tables are checked only on complete
blocks that lie in their stated theorem range. The report explicitly
distinguishes larger parameters whose sufficient transition threshold
lies beyond the finite run. The low-height and elementary Gaussian
identities have separate word-level audits. The proofs, not the finite
checks, establish the assertions for arbitrary indices and parameters;
sampled figure errors are numerical diagnostics, not uniform bounds.

An effective rate of convergence to the elementary Gaussian profile
remains to be determined. A second-order interior expansion would clarify
the visible finite-level differences. The fixed-$p$ results do not give
an expansion uniform in a parameter growing with $n$. Nor do they classify
all seeds producing finite signed concatenations: the theorem concerns
the explicit repeated-pair family. Its common geometric frequency limit
and different finite upper tails show why a histogram alone cannot
recover the underlying ordered dynamics.

\clearpage
\appendix
\section{The three-arch member \texorpdfstring{$Q_2$}{Q2}}\label{app:three-arch}
This appendix records the explicit formulas and finite examples for
$p=2$. Only here we abbreviate
$Q_2,F_2,G_2,U_{2,k},L_{2,k},a_{2,k},d_{2,k},s_{2,k}$
by $Q,F,G,U_k,L_k,a_k,d_k,s_k$, and put $p_k=L_k/2$.
All general inputs and their proofs appear in the main text.

\subsection{Initial block and three-word organization}
\begin{equation}\label{eq:Q}
\begin{split}
 (Q(1),\ldots,Q(7))&=(1,1,3,3,3,3,5),\\
 Q(n)&=Q(n-Q(n-1))+Q(n-Q(n-2))-(-1)^n\qquad(n\ge8).
\end{split}
\end{equation}
The first terms are
\[
1,1,3,3,3,3,5,5,7,5,7,7,7,9,9,9,11,11,11,13,\ldots.
\]
Theorem~\ref{thm:family-global} and Corollary~\ref{cor:family-envelope} give
\begin{equation}\label{eq:basic}
 1\le Q(n)\le n,\qquad Q(n+2)-Q(n)\in\{0,2\},
\end{equation}
\begin{equation}\label{eq:growth}
 Q(n)=n/2+O(n/\sqrt{\log n}),
\end{equation}
with a sharp error order. Every positive odd value is assumed finitely often.
Here the skeleton is
\begin{equation}\label{eq:skeleton}
 a_k=\frac{7\cdot2^k-3-(-1)^k}{6},\qquad
 d_k=\left\lfloor\frac{k-1}{2}\right\rfloor,\qquad s_k=2a_k-d_k.
\end{equation}
The first boundaries are $4,8,17,35,72,146,295$.
The specialized core formulas \eqref{eq:family-core-odd}--\eqref{eq:family-core-even}
and reconstruction \eqref{eq:family-lift} will be used throughout.

\begin{figure}[htbp]
\centering
\includegraphics[width=.98\textwidth]{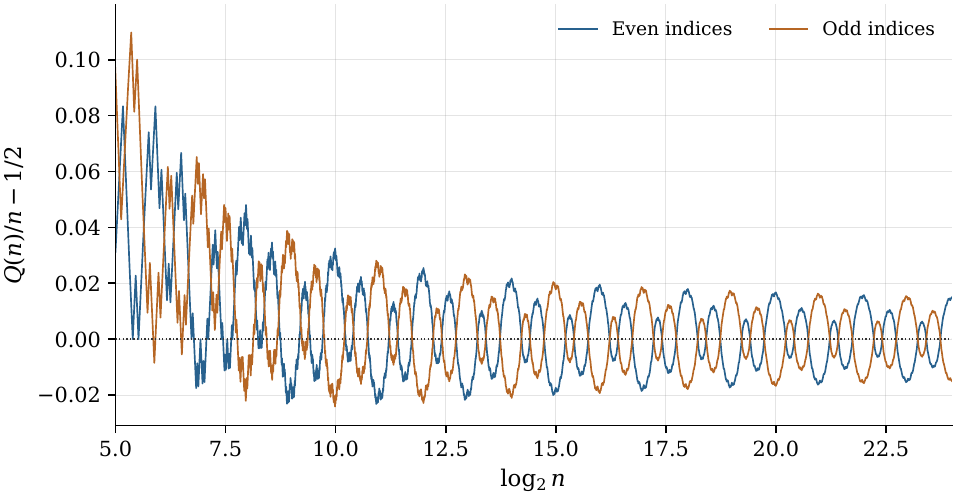}
\caption{The directly computed ratio errors, separated by index parity,
through $n=2^{24}$. The horizontal axis is logarithmic. Several distinct
excursions occur during each unit interval of this axis; the alternating
sign convention is made explicit in Theorem~\ref{thm:standard}.}
\label{fig:overview}
\end{figure}

\begin{theorem}[Exact three-word decomposition]\label{thm:factor}
For every $k\ge2$,
\begin{equation}\label{eq:factor}
U_k=A_k\,\overline{A_{k+1}}\,A_{k+2}.
\end{equation}
Consequently $h_{U_k}$ vanishes exactly at the four integer positions
\begin{equation}\label{eq:word-zeros}
0,\qquad \ell_k,\qquad \ell_k+\ell_{k+1}=2^k,\qquad L_k.
\end{equation}
It is positive, negative, and positive in the three open subintervals.
The three lengths have asymptotic ratio $1:2:4$.
\end{theorem}
\begin{proof}
Specialize \eqref{eq:family-factor} at $p=2$ and use the primitivity
of each $A_k$ and $\ell_k+\ell_{k+1}=2^k$.
\end{proof}
\begin{figure}[htbp]
\centering
\includegraphics[width=.95\textwidth]{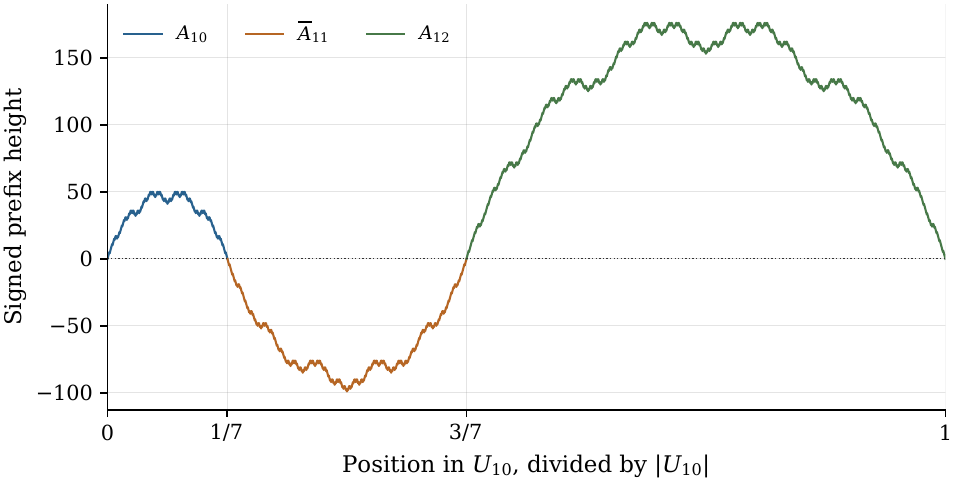}
\caption{The three factors of $U_{10}$ displayed separately. Their
normalized joining positions tend to $1/7$ and $3/7$. Each elementary
excursion is exactly symmetric, although the whole signed concatenation
is not symmetric.}
\label{fig:factor}
\end{figure}

\subsection{The explicit four-piece profile and envelope}\label{sec:profile}

For $1\le t\le2$, define
\begin{equation}\label{eq:profile}
\mathcal P(t)=
\begin{cases}
\displaystyle\frac1{6t}J(3t-5/2),&1\le t\le7/6,\\[4pt]
\displaystyle-\frac1{12t}J(6t-7),&7/6\le t\le4/3,\\[4pt]
\displaystyle\frac1{6t}J(3t-4),&4/3\le t\le5/3,\\[4pt]
\displaystyle-\frac1{3t}J((3t-5)/2),&5/3\le t\le2.
\end{cases}
\end{equation}
The formulas agree at each joining point. 
We use the parity sampling convention from Section~\ref{sec:family-profiles}.

\begin{theorem}[Standard dyadic profile]\label{thm:standard}
For each $r\in\{0,1\}$,
\begin{equation}\label{eq:standard}
\sup_{1\le t\le2}
\left|(-1)^{j+r}\sqrt{j}
\left(\frac{Q(N_{j,r}(t))}{N_{j,r}(t)}-\frac12\right)
-\mathcal P(t)\right|\longrightarrow0.
\end{equation}
The limiting profile has precisely the three zeros
\begin{equation}\label{eq:three-transitions}
7/6,\qquad4/3,\qquad5/3,
\end{equation}
and signs $+,-,+,-$ on the four complementary intervals.
For each parity, all finite-level zeros and sign changes are confined
to neighborhoods of these three positions whose normalized widths tend
to zero. Each neighborhood contains a sign change for all sufficiently
large levels.
\end{theorem}
\begin{proof}
The core formulas and the unit-step bounds imply, uniformly in a macro
block and with $n=2m$ or $2m-1$,
\begin{equation}\label{eq:error-height}
\frac{Q(n)}n-\frac12
=\frac{(-1)^{k+r+1}}{2m}\,h_{U_k}(m-s_k)
+O(k/2^k),
\end{equation}
where $r=n\bmod2$. On even macro levels the shifted bit in $G$ adds
only a bounded error; using $n=2m-1$ rather than $2m$ does the same.

The natural macro boundaries satisfy
$2s_k=(14/3)2^k+O(k)$.
The factorization~\eqref{eq:factor} gives the following four elementary
arcs intersecting a standard block. The endpoints in this table are
limiting positions in the coordinate $t=n/2^j$.
\begin{center}
\begin{tabular}{@{}ccccc@{}}\toprule
Observed interval & Macro level & Word & Full arc interval & Sign\\\midrule
$[1,7/6]$ & $j-3$ & $A_{j-1}$ & $[5/6,7/6]$ & $+$\\
$[7/6,4/3]$ & $j-2$ & $A_{j-2}$ & $[7/6,4/3]$ & $-$\\
$[4/3,5/3]$ & $j-2$ & $\overline{A}_{j-1}$ & $[4/3,5/3]$ & $+$\\
$[5/3,2]$ & $j-2$ & $A_j$ & $[5/3,7/3]$ & $-$\\\bottomrule
\end{tabular}
\end{center}
Signs include the compensation $(-1)^{j+r}$.
An elementary arc of limiting horizontal width $w$, starting at $b$,
contributes $\pm wJ((t-b)/w)/(2t)$ by Theorem~\ref{thm:gaussian} and
\eqref{eq:error-height}. This is precisely \eqref{eq:profile}.
All actual boundaries differ from their tabulated positions by
$O(j2^{-j})$. Uniform elementary convergence, $J(0)=J(1)=0$, and the
$O(j/2^j)$ reconstruction error give uniform convergence through those
boundary neighborhoods as well.

Strict positivity of $J$ in $(0,1)$ proves the sign and zero statements.
Uniform convergence excludes sign changes on compact subintervals away
from the three limiting zeros. Opposite signs on either side of each
zero guarantee at least one change nearby. This does not assert an
exactly three-crossing count at a finite level.
\end{proof}

\begin{remark}[Three transitions, not a factor-three renormalization]
The evolution doubles the elementary lengths at each level, up to
bounded corrections. The three transition points arise from concatenating
three different dyadic generations. They do not describe a triadic
scaling law.
\end{remark}

\begin{figure}[htbp]
\centering
\includegraphics[width=.99\textwidth]{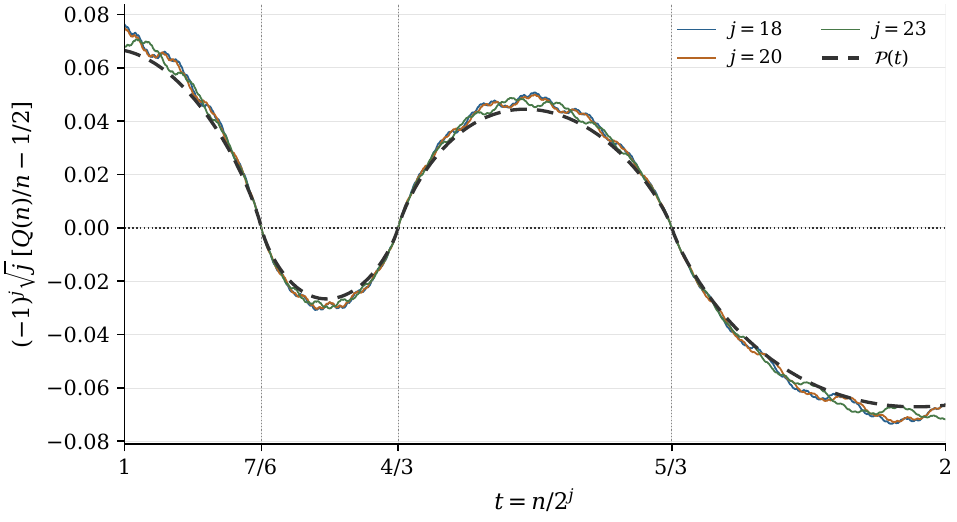}
\caption{Even-index profiles on three standard dyadic blocks and the
explicit limit \eqref{eq:profile}. Only the prescribed sign and
$\sqrt j$ normalization are applied. In particular, the three vertical
reference lines have not been fitted to the data.}
\label{fig:standard}
\end{figure}

\begin{theorem}[Sharp relative-error envelope]\label{thm:envelope}
Let $b_*>0$ be the unique positive solution of
\begin{equation}\label{eq:root}
b_*\bigl(5/2+\Phi(-b_*)\bigr)=\varphi(b_*).
\end{equation}
Then
\begin{equation}\label{eq:envelope}
\limsup_{n\to\infty}\sqrt{\log_2n}
\left|\frac{Q(n)}n-\frac12\right|
=\Lambda_M:=\frac{b_*}{2}
=0.067088003576\ldots.
\end{equation}
The corresponding signed limsup and liminf are $\Lambda_M$ and
$-\Lambda_M$.
\end{theorem}
\begin{proof}
In one complete macro block, the start position of each elementary arc
divided by its own length tends respectively to $7$, $4$, and $5/2$.
Equation~\eqref{eq:error-height} and Theorem~\ref{thm:gaussian} therefore
reduce the envelope to the largest of
\[
\max_{0\le x\le1}\frac{J(x)}{2(\alpha+x)},
\qquad \alpha\in\{7,4,5/2\}.
\]
The expression decreases strictly with $\alpha$, so the third arc wins.
The derivative numerator $(\alpha+x)J'(x)-J(x)$ decreases strictly,
with derivative $(\alpha+x)J''(x)<0$, and crosses zero once.
Writing $x_* =\Phi(-b_*)$ gives \eqref{eq:root}; at that point the
maximum equals $b_*/2$.
The left side minus the right side of \eqref{eq:root} has derivative
$5/2+\Phi(-b)>0$, is negative at zero, and eventually positive.
This proves uniqueness. Opposite parity signs give both signed extremes.
Numerically,
\[
x_*=0.446631698928\ldots,\qquad
5/3+(2/3)x_*=1.964421132618\ldots.
\]
The second number is the maximizing location on the last arc in a
standard block. Finally $\log_2 n=j+O(1)$ on that block, so its square
root may replace $\sqrt j$ in the limiting envelope.
\end{proof}

Equations~\eqref{eq:standard} and \eqref{eq:envelope} give the explicit
three-arch specialization of the general profile and envelope theorems.
In particular, the error cannot be improved to $o(n/\sqrt{\log n})$.

\subsection{All finite zeros, tangencies, and crossings}\label{sec:boundaries}

\begin{theorem}[Every finite zero and every odd-parity crossing]
\label{thm:finite-zeros}
For a standard block $2^j\le n<2^{j+1}$, $j\ge10$, put
\[
 c_{j,1}=\frac76\,2^j,\qquad
 c_{j,2}=\frac43\,2^j,\qquad
 c_{j,3}=\frac53\,2^j.
\]
All solutions of $Q(n)=n/2$ in this block are the even integers
$c_{j,i}+\delta$, with offsets in the following table:
\begin{center}
\begin{tabular}{@{}ccc@{}}\toprule
Center & $j$ even: zero offsets & $j$ odd: zero offsets\\\midrule
$c_{j,1}$ & $\{10/3\}$ & $\{-2j-4/3,\,-2j+8/3\}$\\
$c_{j,2}$ & $\{-2j-10/3\}$ & $\{-2/3,\,10/3\}$\\
$c_{j,3}$ & $\{-2/3,\,10/3\}$ & $\{-2j-4/3,\,-2j+8/3\}$\\\bottomrule
\end{tabular}
\end{center}
In particular, there are exactly four zeros when $j$ is even and six
when $j$ is odd. No odd index can be an exact zero.

On the odd subsequence, there are exactly three sign changes of
$Q(n)-n/2$. They occur between $n$ and $n+2$, where the left endpoint
is $n=c_{j,i}+\eta_{j,i}$ and
\begin{center}
\begin{tabular}{@{}ccc@{}}\toprule
Center & $j$ even: $\eta_{j,i}$ & $j$ odd: $\eta_{j,i}$\\\midrule
$c_{j,1}$ & $-2j-11/3$ & $11/3$\\
$c_{j,2}$ & $11/3$ & $-2j-17/3$\\
$c_{j,3}$ & $-2j-11/3$ & $11/3$\\\bottomrule
\end{tabular}
\end{center}
Every displayed index is an integer of its stated parity.
On the even subsequence there are also exactly three sign changes
when zeros are deleted before counting changes; some of its additional
zeros are tangencies rather than sign changes.
\end{theorem}
\begin{proof}
Write $m=s_k+t$, $d=d_k$, and $h=h_{U_k}$. The core reconstruction
gives the exact identities
\begin{align}
 Q(2m)-m&=
 \begin{cases}d+1+h(t),&k\text{ odd},\\
 d+2-h(t+1),&k\text{ even},\end{cases}
 \label{eq:exact-even-error}\\
 Q(2m-1)-(2m-1)/2&=
 \begin{cases}d+3/2-h(t),&k\text{ odd},\\
 d+3/2+h(t),&k\text{ even}.
 \end{cases}\label{eq:exact-odd-error}
\end{align}
For the even-level formula use
$2U_k[t]=1+h(t)-h(t+1)$. We assign a common endpoint to the block
on its right; the formulas agree there. An odd index cannot be a zero,
since $Q(n)$ is integral and $n/2$ is then a half-integer.

For an elementary word with initial run length $\lambda$, the
low-height lemma gives the following complete lists:
\[
\begin{array}{c|c}
\text{height}&\text{positions}\\\hline
\lambda&\lambda,\ \lambda+2,\ \ell-\lambda-2,\ \ell-\lambda\\
\lambda+1&\lambda+3,\ \ell-\lambda-3.
\end{array}
\]
It also determines all crossings of the half-integer thresholds used
below. The upward crossings of $\lambda+1/2$ and $\lambda+3/2$
have left endpoints $\lambda+2$ and $\lambda+3$, respectively;
their downward counterparts have left endpoints
$\ell-\lambda-3$ and $\ell-\lambda-4$. There are no other
crossings of these thresholds: the middle heights are at least
$\lambda+3$.

Here is the full coordinate calculation. Set
$r=\lfloor(j-2)/2\rfloor$. The four elementary factors meeting the
standard block have core-coordinate beginnings
\[
\begin{array}{c|c|c|c}
\text{factor}&\text{macro level}&\text{word}&\text{beginning }B\\\hline
\mathrm I&j-3&A_{j-1}&s_{j-3}+2^{j-3}\\
\mathrm {II}&j-2&A_{j-2}&s_{j-2}\\
\mathrm {III}&j-2&\overline{A}_{j-1}&s_{j-2}+\ell_{j-2}\\
\mathrm {IV}&j-2&A_j&s_{j-2}+2^{j-2}.
\end{array}
\]
In the next two tables, $(\mathrm X,q)$ denotes the core index
$m=B_{\mathrm X}+q$. The entries in the first table are exactly the
even-index zeros, with $n=2m$:
\[
\begin{array}{c|c|c}
\text{transition}&j\text{ even}&j\text{ odd}\\\hline
1&(\mathrm {II},r+2)&
 (\mathrm I,\ell_{j-1}-r-4),\ (\mathrm I,\ell_{j-1}-r-2)\\
2&(\mathrm {II},\ell_{j-2}-r-4)&
 (\mathrm {III},r+1),\ (\mathrm {III},r+3)\\
3&(\mathrm {IV},r),\ (\mathrm {IV},r+2)&
 (\mathrm {III},\ell_{j-1}-r-3),\ (\mathrm {III},\ell_{j-1}-r-1).
\end{array}
\]
For odd indices, the entries below give the left member $n=2m-1$
of the pair $n,n+2$ across which the sign changes:
\[
\begin{array}{c|c|c}
\text{transition}&j\text{ even}&j\text{ odd}\\\hline
1&(\mathrm I,\ell_{j-1}-r-3)&(\mathrm {II},r+3)\\
2&(\mathrm {III},r+2)&(\mathrm {II},\ell_{j-2}-r-4)\\
3&(\mathrm {III},\ell_{j-1}-r-3)&(\mathrm {IV},r+3).
\end{array}
\]
To obtain these tables, on an odd macro level apply the even-error
identity only to the negative middle excursion, at height $d+1$.
On an even macro level apply it to the positive first and third
excursions, at height $d+2$, remembering the shift from $t$ to $t+1$.
For the odd-error identity, positive excursions on odd macro levels
and negative excursions on even macro levels cross height $d+3/2$.
The complete low-height lists above then give every displayed entry.

Substituting the explicit formulas for $s_k$ and $\ell_k$ converts
these core coordinates into the two offset tables in the statement.
For example, for even $j$ the first even zero is
\[
 2(s_{j-2}+r+2)=\frac76\,2^j+\frac{10}{3},
\]
while the first odd crossing has left endpoint
\[
 2(s_{j-3}+2^{j-3}+\ell_{j-1}-r-3)-1
 =\frac76\,2^j-2j-\frac{11}{3}.
\]
The exact same substitutions in the remaining explicitly listed
coordinates give the other offsets without an asymptotic argument.

All elementary levels used here are at least eight when $j\ge10$.
The initial point of the standard block lies strictly inside factor
$\mathrm I$, and its final point strictly inside factor $\mathrm {IV}$,
beyond the short end intervals adjacent to the three intervening
joins. These assertions follow directly from the displayed beginnings
and lengths; the distance to the nearest excluded end is of order
$2^j$, whereas each relevant end interval has length at most $j/2+4$.
For an explicit uniform separation, $2^j/6>4j+16$ holds at $j=10$
and persists on increasing $j$. Thus all the listed crossings lie
in the standard block, the three groups are disjoint, and no omitted
outer endpoint contributes a zero or crossing.

Finally, the local pattern $g_\lambda$ distinguishes tangencies from
genuine changes on the even subsequence. At height $\lambda$, the
first hit $\lambda$ is followed by a return to $\lambda-1$, whereas
the second hit $\lambda+2$ is followed by heights greater than
$\lambda$; reflection gives the other end. At height $\lambda+1$
there is a single transverse hit at each end. Consequently deleting
the zeros leaves exactly one sign change near each of the three
centers. The low-height lemma excludes any further one.
\end{proof}

\paragraph{A worked example: the block $1024\le n<2048$.}
For $j=10$, the three limiting centers are
\[
 c_{10,1}=\frac{3584}{3},\qquad
 c_{10,2}=\frac{4096}{3},\qquad
 c_{10,3}=\frac{5120}{3}.
\]
The exact word coordinates retain the small offsets that these centers
do not record. Here $r=4$, $d_7=d_8=3$, and
\[
 s_7=295,\quad s_8=593,\qquad
 \ell_8=86,\quad\ell_9=170,\quad\ell_{10}=342.
\]
Thus the four factors in the proof of Theorem~\ref{thm:finite-zeros}
start at core indices
\[
 B_{\mathrm I}=423,\quad B_{\mathrm {II}}=593,\quad
 B_{\mathrm {III}}=679,\quad B_{\mathrm {IV}}=849.
\]
They are respectively $A_9$, $A_8$, $\overline A_9$, and $A_{10}$.
Write $E(n)=Q(n)-n/2$.

On factors $\mathrm {II}$ and $\mathrm {IV}$, the even-error formula is
$E(2m)=5-h_{U_8}(m-s_8+1)$. In factor $\mathrm {II}$ it vanishes at
local coordinates $q=6,78$: the two prefix positions are $q+1=7,79$,
and
\[
 h_{A_8}(7)=h_{A_8}(86-7)=g_4(7)=5.
\]
In factor $\mathrm {IV}$ it vanishes at $q=4,6$, since
\[
 h_{A_{10}}(5)=g_5(5)=5,
 \qquad h_{A_{10}}(7)=g_5(7)=5.
\]
The other end of factor $\mathrm {IV}$ is beyond this standard block.
Factor $\mathrm I$ has $E(2m)=4+h_{A_9}(q)>0$, and factor
$\mathrm {III}$ has $E(2m)=5+h_{A_9}(q+1)>0$.
The low-height lemma therefore leaves exactly the following four zeros.
The adjacent values distinguish a tangency from a genuine crossing:
\begin{center}
\small
\begin{tabular}{@{}ccccc@{}}\toprule
$(\mathrm X,q)$ & $n$ & $(Q(n-2),Q(n),Q(n+2))$
 & $(E(n-2),E(n),E(n+2))$ & Type\\\midrule
$(\mathrm {II},6)$ & $1198$ & $(599,599,599)$ & $(1,0,-1)$ & crossing\\
$(\mathrm {II},78)$ & $1342$ & $(669,671,673)$ & $(-1,0,1)$ & crossing\\
$(\mathrm {IV},4)$ & $1706$ & $(853,853,855)$ & $(1,0,1)$ & tangency\\
$(\mathrm {IV},6)$ & $1710$ & $(855,855,855)$ & $(1,0,-1)$ & crossing\\\bottomrule
\end{tabular}
\end{center}
For example, the first zero is $2(593+6)=1198=c_{10,1}+10/3$.
The remaining offsets are $-70/3$ at $c_{10,2}$ and $-2/3,10/3$
at $c_{10,3}$. In particular, the two zeros near the third center
produce only one sign change after zeros are deleted.

For odd indices in factor $\mathrm I$, the exact error is
$9/2-h_{A_9}(q)$. The reflected low heights at $q=163,164$ are $5,4$,
so the crossing occurs between $2(423+163)-1=1171$ and $1173$.
In factor $\mathrm {III}$ the exact error is again
$9/2-h_{A_9}(q)$, now because the factor is complemented and the
macro level is even. The pairs of local coordinates $(6,7)$ and
$(163,164)$ have heights $(4,5)$ and $(5,4)$, respectively. Thus
\begin{center}
\begin{tabular}{@{}cccc@{}}\toprule
$(\mathrm X,q)$ & $(n,n+2)$ & $(Q(n),Q(n+2))$ & $(E(n),E(n+2))$\\\midrule
$(\mathrm I,163)$ & $(1171,1173)$ & $(585,587)$ & $(-1/2,1/2)$\\
$(\mathrm {III},6)$ & $(1369,1371)$ & $(685,685)$ & $(1/2,-1/2)$\\
$(\mathrm {III},163)$ & $(1683,1685)$ & $(841,843)$ & $(-1/2,1/2)$\\\bottomrule
\end{tabular}
\end{center}
These are the three odd-parity crossings. The left-endpoint offsets
are $-71/3$, $11/3$, and $-71/3$, respectively, as predicted by the
theorem. The archive also checks every value through
$n=2048$ by direct recursion and, independently, by word reconstruction;
both calculations give the values in these tables.

\begin{corollary}[An optimal logarithmic transition window]
\label{cor:boundary-rate}
All zeros and parity-separated sign transitions in the $j$-th standard
block lie in the union of the intervals
\begin{equation}\label{eq:boundary-window}
 |n-c_{j,i}|\le2j+8,\qquad i=1,2,3,\quad j\ge10.
\end{equation}
Thus their distance from the limiting positions in the coordinate
$n/2^j$ is $O(j2^{-j})=O((\log n)/n)$.
This order cannot uniformly be replaced by $o(j2^{-j})$, as the offsets
$-2j+O(1)$ in Theorem~\ref{thm:finite-zeros} show.
\end{corollary}

\subsection{Natural and pure dyadic frequency blocks}\label{sec:frequencies}
Write $\nu(v)=\#\{n\ge1:Q(n)=v\}$; the subscript in
$\nu_2(b)$ below instead denotes the two-adic valuation.
Put $\mathcal V_k=\{2u+1:a_k\le u<a_{k+1}\}$ and
$\mathcal F_k(z)=\sum_{u=a_k}^{a_{k+1}-1}z^{\nu(2u+1)}$.

\begin{theorem}[Complete-block multiplicity recursion]\label{thm:frequency}
Every coefficient of $\mathcal F_k$ is determined by
\begin{equation}\label{eq:frequency-rec}
\mathcal F_1(z)=z^3+z^4,\qquad
\mathcal F_{k+1}(z)=z\mathcal F_k(z)+(p_k-(-1)^k)z^3
\quad(k\ge1).
\end{equation}
Consequently
\begin{equation}\label{eq:pgf-limit}
\frac{\mathcal F_k(z)}{p_k}\longrightarrow
G(z)=\frac{z^3}{2-z}
\end{equation}
in coefficient $\ell^1$ norm with error $O(k2^{-k})$.
The limiting distribution has mean $4$ and variance $2$; the complete
block means and variances converge to these values.
For every $k\ge1$,
\begin{align}
\sum_{v\in\mathcal V_k}\nu(v)&=4p_k-\frac{3+(-1)^k}{2},\label{eq:frequency-mass}\\
\max_{v\in\mathcal V_k}\nu(v)&=k+3.\label{eq:frequency-max}
\end{align}
In particular, the multiplicities are unbounded but are $O(\log v)$.
\end{theorem}
\begin{proof}
Specialize Theorem~\ref{thm:family-natural-frequency} at $p=2$.
Here the explicit first step also satisfies the same polynomial recursion:
$\mathcal F_1=z^3+z^4$ and $\mathcal F_2=3z^3+z^4+z^5$.
All stated masses, degrees, and convergence estimates follow.
\end{proof}
The ordered natural-block statement becomes
\begin{equation}\label{eq:ordered-natural-frequencies}
 \bigl(\nu(2u+1)\bigr)_{u=a_k}^{a_{k+1}-1}
 =\mathbf V_k\mathbf V_{k+1}\mathbf V_{k+2}.
\end{equation}
There is no repeated-pair join and hence no tilde correction.

\begin{theorem}[Pure dyadic value blocks]\label{thm:power-blocks}
For every $j\ge1$, as multisets,
\begin{equation}\label{eq:power-law}
 \{\nu(2u+1):2^j\le u<2^{j+1}\}
 =\{3+\nu_2(b):1\le b\le2^j\},
\end{equation}
where $\nu_2$ is the exponent of two dividing a positive integer.
Equivalently, the exact polynomial is
\begin{equation}\label{eq:power-polynomial}
 D_j(z):=\sum_{u=2^j}^{2^{j+1}-1}z^{\nu(2u+1)}
 =\sum_{h=0}^{j-1}2^{j-h-1}z^{h+3}+z^{j+3}.
\end{equation}
More precisely, the \emph{ordered} frequency word on this block is,
for $j\ge2$,
\begin{equation}\label{eq:ordered-power-word}
 \operatorname{rev}(\mathbf R_{j+1})\,
 \mathbf V_j\,\mathbf V_{j+1}\,(j+3)\,\mathbf R_{j+2}.
\end{equation}
\end{theorem}
\begin{proof}
For $j\ge3$ the lower cut $2^j$ lies in the last elementary frequency
factor of the preceding macro block. That factor is $\mathbf V_{j+1}$
and starts at
\[
 b_j=a_{j-1}+q_{j-1}+q_j
     =\frac{5\cdot2^j-3+(-1)^j}{6}.
\]
Since
\begin{equation}\label{eq:dyadic-center-cut}
 2^j=b_j+\frac{q_{j+1}+1}{2},
\end{equation}
its retained suffix is exactly $\operatorname{rev}(\mathbf R_{j+1})$.
The next macro block contributes the whole words $\mathbf V_j$ and
$\mathbf V_{j+1}$. Applying \eqref{eq:dyadic-center-cut} at the upper
endpoint shows that the final retained prefix is
$(j+3)\mathbf R_{j+2}$. This proves
\eqref{eq:ordered-power-word}; its $j=2$ instance is the direct word
$(3,4,5,3)$.

The polynomial of $\mathbf R_k$ is
$(\mathcal B_k(z)-z^{k+1})/2$. Hence, for $j\ge2$,
\begin{equation}\label{eq:D-from-E}
 2D_j=2\mathcal B_j+3\mathcal B_{j+1}+\mathcal B_{j+2}-z^{j+2}+z^{j+3}.
\end{equation}
Write $e=(-1)^j$. From $q_{j+1}=2q_j-e$ and
$q_{j+2}=4q_j-e$, the non-shifted terms when advancing the three elementary polynomials
sum to
\[
 (q_j-e)+\tfrac32(q_{j+1}+e)+\tfrac12(q_{j+2}-e)
 =6q_j-2e=2^j.
\]
The two displayed monomial corrections in \eqref{eq:D-from-E} simply
shift by $z$. We therefore obtain
\begin{equation}\label{eq:D-recursion}
 D_{j+1}(z)=zD_j(z)+2^jz^3.
\end{equation}
The direct base $D_1=z^3+z^4$ also gives
$D_2=2z^3+z^4+z^5$ by the same rule.
Solving this scalar recursion proves \eqref{eq:power-polynomial} and
therefore \eqref{eq:power-law}. No block-mass or hitting-time assumption
is used in the proof.
\end{proof}

The proof of Theorem~\ref{thm:power-blocks} depends on the locations
of the cuts through the palindromic tails. The natural-block polynomial
alone does not determine these cuts. The valuation law describes the
multiset; it does not assert the pointwise formula
$\nu(2u+1)=3+\nu_2(u)$.

\Needspace{14\baselineskip}
\begin{corollary}[Exact moments, convergence rate, and cumulative count]
\label{cor:power-statistics}
Choose a value uniformly from
$\{2u+1:2^j\le u<2^{j+1}\}$, $j\ge1$, and let $X_j$ be its total
multiplicity in $Q$. Then
\begin{align}
 \mathbb E X_j&=4-2^{-j},\label{eq:power-mean}\\
 \operatorname{Var}(X_j)&=2-(2j+1)2^{-j}-4^{-j}.
 \label{eq:power-variance}
\end{align}
Its law has coefficient $\ell^1$ distance exactly $2^{-j}$ from
$G(z)=z^3/(2-z)$, and total-variation distance exactly $2^{-j-1}$.
Also
\begin{equation}\label{eq:threshold-count}
 \#\{n\ge1:Q(n)<2^{j+1}\}=2^{j+2}-j-1.
\end{equation}
This is a count over the whole infinite sequence, not a claim that the
counted indices form an initial interval.
\end{corollary}
\begin{proof}
The distribution of $X_j-3$ is that of a geometric random variable
with probabilities $2^{-h-1}$ at $h\ge0$, truncated by replacing every
value at least $j$ by $j$. Its first two moments, or differentiation
of \eqref{eq:power-polynomial}, give
\eqref{eq:power-mean}--\eqref{eq:power-variance}.
At multiplicity $j+3$ there is an excess probability $2^{-j-1}$;
the removed tail has that same mass. This proves the two distances.
The exact block mass is $4\cdot2^j-1$.
Finally $\nu(1)=2$ and $\nu(3)=4$; summing those two exceptions and
all blocks below $2^j$ proves \eqref{eq:threshold-count}.
\end{proof}

\begin{proposition}[Locations of the two largest multiplicities]
\label{prop:power-extrema}
On $2^j\le u<2^{j+1}$, $j\ge1$, the largest multiplicity $j+3$ occurs
at the unique value index
\begin{equation}\label{eq:power-max-location}
 u=\left\lfloor\frac{5\cdot2^j}{3}\right\rfloor.
\end{equation}
The second largest multiplicity $j+2$ occurs uniquely at
$u=\lfloor4\cdot2^j/3\rfloor$.
\end{proposition}
\begin{proof}
By Lemma~\ref{lem:frequency-words}, the only entry $k+1$ in
$\mathbf V_k$ is its first, and none of its entries equals $k$.
Thus in \eqref{eq:ordered-power-word} the two largest entries are
the first entry of $\mathbf V_{j+1}$ and the displayed $j+3$.
Their positions in the macro block are
\[
 a_j+q_j=\left\lfloor\frac{4\cdot2^j}{3}\right\rfloor,
 \qquad
 a_j+q_j+q_{j+1}=\left\lfloor\frac{5\cdot2^j}{3}\right\rfloor.
\]
This proves the claim for $j\ge2$; the block $j=1$ is immediate.
\end{proof}

\begin{figure}[htbp]
\centering
\includegraphics[width=.96\textwidth]{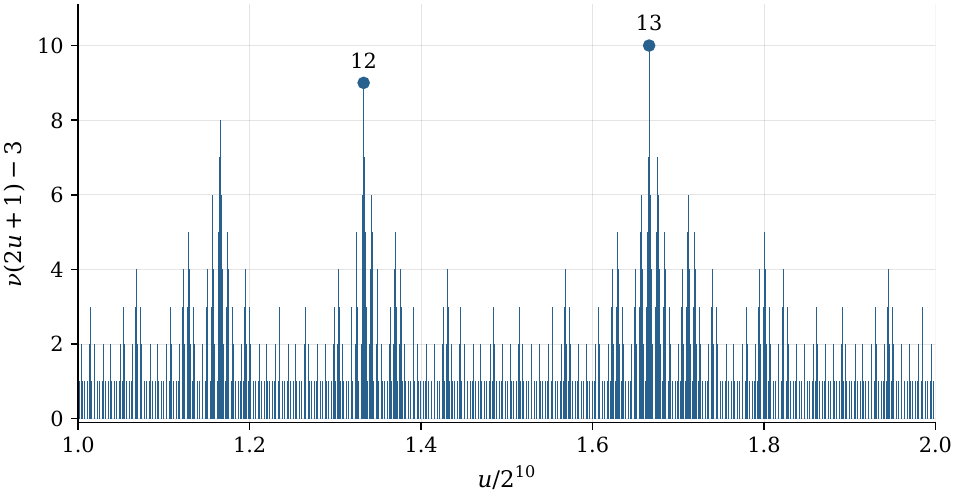}
\caption{Ordered multiplicities in a pure dyadic value-index block.
The histogram is exactly the shifted valuation law of
Theorem~\ref{thm:power-blocks}, but the order is not the ruler order.
The two marked locations are the exact maximizers specified in
Proposition~\ref{prop:power-extrema}. Here the horizontal variable is a
\emph{value index} $u$, not the sequence index $n$, and $j=10$.
The labels at the two marked points give their multiplicities.}
\label{fig:power-frequencies}
\end{figure}

\subsection{Exact summation in the three-arch case}\label{sec:summatory}

The large, opposite-parity oscillations do not average to zero exactly.
The symmetric core mode determines a smaller systematic drift. It also
allows all complete macro blocks to be summed without evaluating their
individual heights.

Let
\[
 S(N)=\sum_{n=1}^N Q(n),\qquad
 \mathcal E(N)=S(N)-\frac{N(N+1)}4.
\]
The subtraction in $\mathcal E$ is the exact sum of the comparison
terms $n/2$.

\begin{theorem}[An exact telescoping formula]\label{thm:summatory-exact}
Let $M\ge4$, and let $K$ be the unique integer for which
$s_K\le M<s_{K+1}$. Put $T=M-s_K+1$. Then
\begin{equation}\label{eq:exact-sum}
 \mathcal E(2M)
 =\left(K+\frac32\right)(M+1)
 -\sum_{k=1}^K s_k-\frac52-d_K
 -\mathbf1_{\{K\text{ even}\}}h_{U_K}(T).
\end{equation}
The sum of the boundary positions is itself explicit:
\begin{equation}\label{eq:sum-boundaries}
 \sum_{k=1}^K s_k
 =\frac73(2^{K+1}-2)-K-\frac{(-1)^K-1}{6}
 -\left\lfloor\frac{(K-1)^2}{4}\right\rfloor.
\end{equation}
Thus only one height of the final, incomplete word is needed beyond
powers, floors, and elementary arithmetic. When $K$ is odd, even that
height term is absent.
\end{theorem}
\begin{proof}
For $m=s_k+t$ the exact symmetric-mode formula gives the contribution
of the pair $2m-1,2m$ to $\mathcal E$ as
\[
 2\bigl(F(m)+G(m)-m\bigr)+\frac52
 =\begin{cases}
 k+3/2,&k\text{ odd},\\
 k+3/2-(h_{U_k}(t+1)-h_{U_k}(t)),&k\text{ even}.
 \end{cases}
\]
The first three pairs contribute $7/2$ in total.
Assign each macro block the half-open index interval
$[s_k,s_{k+1})$. A completed odd block therefore contributes
$(k+3/2)(s_{k+1}-s_k)$.
For an even block, its length is $L_k-1$ and
$h_{U_k}(L_k-1)=1$, because its final bit is one. Its contribution is
$(k+3/2)(s_{k+1}-s_k)-1$.
The unfinished block contributes $(K+3/2)T$ minus the height term
in \eqref{eq:exact-sum}.
There are $d_K$ completed even blocks. Telescoping the terms involving
$s_{k+1}-s_k$, with $s_1=4$, proves \eqref{eq:exact-sum}.
Finally substitute \eqref{eq:skeleton} and use
\[
 \sum_{k=1}^K(-1)^k=\frac{(-1)^K-1}{2},\qquad
 \sum_{k=1}^K\left\lfloor\frac{k-1}{2}\right\rfloor
 =\left\lfloor\frac{(K-1)^2}{4}\right\rfloor
\]
to obtain \eqref{eq:sum-boundaries}.
\end{proof}

\Needspace{16\baselineskip}
\begin{theorem}[A periodic linear term in the summatory asymptotic]
\label{thm:summatory-asymptotic}
Put $C=14/3$, and let $\Psi$ be the continuous one-periodic function
whose values on $0\le\theta\le1$ are
\begin{equation}\label{eq:summatory-periodic}
 \Psi(\theta)=1-\frac12\log_2 C-\frac\theta2-2^{-\theta}.
\end{equation}
Its two endpoint values agree. As $N\to\infty$ through all integers,
\begin{equation}\label{eq:summatory-asymptotic}
 S(N)=\frac{N^2}{4}+\frac N2\log_2 N
      +N\Psi\left(\left\{\log_2(N/C)\right\}\right)
      +O\!\left(\frac N{\sqrt{\log N}}\right).
\end{equation}
Here $\{x\}=x-\lfloor x\rfloor$ denotes the fractional part.
In particular, the cancellation of the large parity oscillations
leaves a positive bias of order $N\log N$, with leading coefficient
$1/2$ when the logarithm has base two.
\end{theorem}
\begin{proof}
First let $N=2M$. By the three-word factorization and the elementary
height bound, uniformly on the $K$-th macro block,
$|h_{U_K}|=O(2^K/\sqrt K)$.
Equations~\eqref{eq:exact-sum}--\eqref{eq:sum-boundaries} yield
\begin{equation}\label{eq:summatory-intermediate}
 \mathcal E(2M)
 =M(K+3/2)-\frac73\,2^{K+1}
  +O\!\left(K^2+\frac{2^K}{\sqrt K}\right).
\end{equation}
Put $b_k=(7/3)2^k$.  The exact boundary formula is
\[
 s_k=b_k-e_k,\qquad
 e_k=1+\frac{(-1)^k}{3}+\left\lfloor\frac{k-1}{2}\right\rfloor>0.
\]
Let $K_0=\lfloor\log_2(M/(7/3))\rfloor$.
For all sufficiently large $M$, either $K=K_0$ or $K=K_0+1$:
the first inequality follows from $s_{K_0}<b_{K_0}\le M$,
and $s_{K_0+2}>b_{K_0+1}>M$ gives the second.
If $K=K_0+1$, then
\[
 b_K-e_K\le M<b_K,\qquad |M-b_K|\le e_K=O(K).
\]
The difference between the main terms of
\eqref{eq:summatory-intermediate} at $K$ and at $K_0$ is exactly
$M-b_K$.  Replacing $K$ by $K_0$ therefore costs only $O(K)$, uniformly
across the short intervals between the actual and ideal dyadic boundaries.

Set $\theta=\{\log_2(N/C)\}$, so
$K_0=\log_2 N-\log_2 C-\theta$ and
$(7/3)2^{K_0+1}=N2^{-\theta}$.
Substituting these in \eqref{eq:summatory-intermediate} and adding
$N(N+1)/4$ gives \eqref{eq:summatory-asymptotic}.
The formula \eqref{eq:summatory-periodic} has equal values at zero
and one, so its periodic extension is continuous and Lipschitz.
For an odd $N$ add its last term to the formula at $N-1$ and use
$Q(N)=N/2+O(N/\sqrt{\log N})$.
The smooth main terms change by $N/2+O(\log N)$, and the periodic term
by $O(1)$. This proves the same estimate for odd $N$.
\end{proof}

\begin{remark}[The finite-level remainder and macro-level parity]
For even $N=2M$, retain the actual level $K$ and $T=M-s_K+1$ from
Theorem~\ref{thm:summatory-exact}, and put
$\theta=\{\log_2(N/(14/3))\}$.  The same calculation, with the final
height kept instead of bounded, gives the uniform refinement
\begin{equation}\label{eq:summatory-height-remainder}
 \frac{S(N)-N^2/4-(N/2)\log_2N}{N}
 =\Psi(\theta)
 -\mathbf1_{\{K\text{ even}\}}\frac{h_{U_K}(T)}N
 +O\!\left(\frac{(\log N)^2}{N}\right).
\end{equation}
Thus odd macro levels have a remainder $O((\log N)^2)$ before
normalization.  At even macro levels the signed word height produces a
larger, visibly oscillatory discrepancy of order at most
$1/\sqrt{\log N}$ after normalization.  Both tend uniformly to zero;
their markedly different finite-level appearance does not require a
different periodic limiting function.
\end{remark}

\begin{figure}[htbp]
\centering
\includegraphics[width=.77\textwidth]{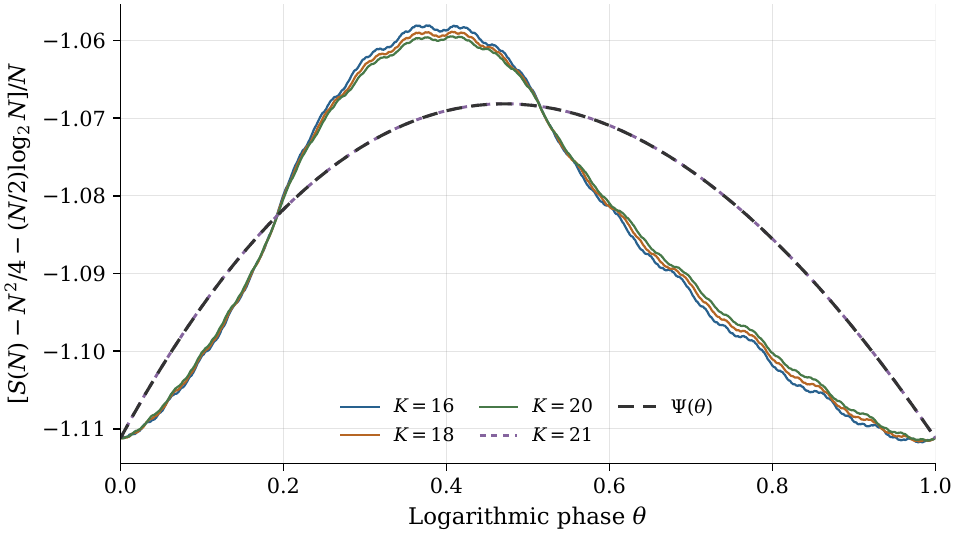}
\caption{The normalized summatory residual and the periodic term $\Psi$.
The phase is sampled at even sequence indices
$N\approx(14/3)2^K2^\theta$. The even scales $K=16,18,20$ show the
signed height correction in \eqref{eq:summatory-height-remainder}.
On the odd scale $K=21$ it is absent except in the short boundary
intervals where the actual macro level differs from $K$; the agreement
with $\Psi$ is consequently much closer. The actual and ideal
boundaries differ by $O(K)$ indices.}
\label{fig:summatory}
\end{figure}
\FloatBarrier

\paragraph{Acknowledgments}
The author thanks Beno\^it Cloitre for his work on the exact-fit
interleaving framework, which provides the combinatorial foundation
for the word constructions used here.

\paragraph{AI-Disclosure}
The author used OpenAI's ChatGPT for help in language editing,
LaTeX restructuring, mathematical explorations, and the preparation of
verification scripts. All mathematical statements, proofs, computations,
references, and the final presentation were independently verified by
the author, who assumes full responsibility for the content.


\begin{thebibliography}{99}

\bibitem{balamohan2007}
B.~Balamohan, A.~Kuznetsov, and S.~Tanny,
\emph{On the behavior of a variant of Hofstadter's $Q$-sequence},
J. Integer Seq. \textbf{10} (2007), Article 07.7.1.
\url{https://cs.uwaterloo.ca/journals/JIS/VOL10/Tanny/tanny3.html}.

\bibitem{bobkov1997}
S.~G. Bobkov,
\emph{An isoperimetric inequality on the discrete cube, and an elementary
proof of the isoperimetric inequality in Gauss space},
Ann. Probab. \textbf{25} (1997), 206--214.
\href{https://doi.org/10.1214/aop/1024404285}{doi:10.1214/aop/1024404285}.

\bibitem{cloitre}
B.~Cloitre,
\emph{The Mantovanelli--Hofstadter Sequence},
arXiv:2604.06237v3 [math.NT] (2026).
\url{https://arxiv.org/abs/2604.06237v3}.

\bibitem{hofstadter1979}
D.~R. Hofstadter,
\emph{G{\"o}del, Escher, Bach: An Eternal Golden Braid},
Basic Books, New York, 1979.

\bibitem{jacksonruskey2006}
B.~Jackson and F.~Ruskey,
\emph{Meta-Fibonacci sequences, binary trees and extremal compact codes},
Electron. J. Combin. \textbf{13} (2006), Research Paper 26.
\href{https://doi.org/10.37236/1052}{doi:10.37236/1052}.

\bibitem{mantovanelli}
M.~Mantovanelli,
\emph{A Dyadic Frequency Law for a Perturbed Hofstadter $Q$-Recursion},
arXiv:2603.16111v3 [math.CO] (2026).
\url{https://arxiv.org/abs/2603.16111v3}.

\bibitem{mantovanellizeta}
M.~Mantovanelli,
\emph{A Two-Variable Zeta Function for a Parity-Perturbed Hofstadter
$Q$-Recursion: The Exceptional $t=-1$ Slice and Gaussian Boundary Layers},
arXiv:2609.02412 [math.NT] (2026).
\url{https://arxiv.org/abs/2609.02412}.

\bibitem{mantovanellidelayed}
M.~Mantovanelli,
\emph{A New Perturbed Hofstadter $Q$-Recursion: Exact Frequencies and
Gaussian Profiles}, ResearchGate preprint (2026).
\href{https://doi.org/10.13140/RG.2.2.36783.70561}
{doi:10.13140/RG.2.2.36783.70561}.
Reproducibility archive:
\href{https://doi.org/10.5281/zenodo.23426341}
{doi:10.5281/zenodo.23426341}.

\bibitem{oeis}
The OEIS Foundation Inc.,
\emph{The On-Line Encyclopedia of Integer Sequences}, entry A394051,
accessed 17 October 2026.
\url{https://oeis.org/A394051}.

\bibitem{oeisA401334}
The OEIS Foundation Inc.,
\emph{The On-Line Encyclopedia of Integer Sequences}, entry A401334
(2026).
\url{https://oeis.org/A401334}.

\bibitem{tanny1992}
S.~M. Tanny,
\emph{A well-behaved cousin of the Hofstadter sequence},
Discrete Math. \textbf{105} (1992), 227--239.
\href{https://doi.org/10.1016/0012-365X(92)90145-6}
{doi:10.1016/0012-365X(92)90145-6}.

\bibitem{zenodoarchive}
M.~Mantovanelli,
\emph{An Odd-Arch Family of Perturbed Hofstadter $Q$-Recursions:
Gaussian Profiles and Exact Frequency Laws (Reproducibility Archive)},
Zenodo (2026).
\href{https://doi.org/10.5281/zenodo.23593788}
{doi:10.5281/zenodo.23593788}.

\end{thebibliography}
\end{document}